\pdfoutput=1
\documentclass{jmlr}
\makeatletter
\def\set@curr@file#1{\def\@curr@file{#1}}
\let\Ginclude@graphics\@org@Ginclude@graphics
\renewcommand\ps@jmlrtps{%
  \let\@mkboth\@gobbletwo
  \def\@oddhead{}\let\@evenhead\@oddhead
  \def\@oddfoot{\hfill\small\rmfamily\thepage\hfill}\let\@evenfoot\@oddfoot}
\makeatother

\title[Error-Corrected Neural Operator for Optimal Control under Uncertainty]{
Error-Corrected Physics-Informed Neural Operator for Optimal Control of Unmanned Aerial Vehicles under Uncertainty}
\usepackage{times}
\usepackage{booktabs}
\usepackage{array}

\author{%
 \Name{Oliver Lundqvist} \Email{oliver.lundqvist@aalto.fi}\\
 \addr Aalto University
 \AND
 \Name{Fabricio Oliveira} \Email{fabol@dtu.dk}\\
 \addr Technical University of Denmark%
}

\begin{document}

\maketitle

\begin{abstract}%
\sloppy We solve open-loop optimal control problems under parametric uncertainty using a physics-informed DeepONet that returns trajectories for all sampled scenarios in one evaluation. Trained only on the residual at generic controls, the operator is reused across targets and horizons. Because a gradient-based optimizer exploits surrogate errors, we correct the prediction by one Newton step on the dynamics, computed for all scenarios and times by one parallel scan. The corrected state has a second-order error bounded by the residual alone, and the corrected cost is honest for every control up to the sampling error. On a fixed-wing UAV benchmark, the corrected operator matches the accuracy of a batched Adam solver and an Ipopt solver at equal sample size and is 8 times faster than the Adam solver on the same GPU at a few hundred scenarios and 15 times faster than Ipopt on a CPU.%
\end{abstract}


\section{Introduction}

Uncertain optimal control refers to an optimal control problem in which the dynamics are subject to parametric uncertainty \citep{ross2015riemann, phelps2016optimal}. The cost function is typically a probabilistic quantity, such as a risk measure. For an expectation, the sample-average approximation replaces it with the empirical mean over sampled scenarios \citep{shapiro2009lectures, phelps2014consistent, melnikov2026convergence}.
The fixed-wing unmanned aerial vehicle (UAV) problem of \citet{lambrianides2020scalable}, with an uncertain start state and uncertain aerodynamic coefficients, requires many scenarios, so every gradient of the sample average integrates the flight equations once per scenario over the horizon, making it costly to compute. In contrast, sigma-point and polynomial-chaos methods need far fewer samples, but their sample count grows with the number of uncertain parameters, which is likely why the earlier study of the same aircraft considered only one uncertain coefficient \citep{shaffer2016unscented}. The question this paper addresses is whether a learned, computationally efficient solution map, evaluated once for all scenarios, can replace the integrator inside the optimizer to reduce the solution time without giving up the guarantees of a time-stepping solver.

Neural operators provide an efficient solution map for the dynamics and act as differentiable surrogates inside the optimizer \citep{wang2021fast}. Computing gradients by automatic differentiation through the neural operator needs caution, since the optimizer exploits surrogate errors \citep{lundqvist2025employing}.  The long horizon of the UAV problem is a hard case, since errors accumulate over time \citep{hairer1993solving, wang2023long, mccabe2023towards} and the optimizer can leave the training distribution \citep{zhu2023reliable}. Remedies include derivative-informed training \citep{oleary2024derivative, luo2025efficient, yao2025derivative}, architectural designs \citep{kurutach2018model, allen2022physical, mandl2026physics}, residual penalization \citep{lundqvist2025employing}, trust regions on certified surrogates \citep{kouri2013trust, qian2017certified}, retraining near the iterate \citep{lye2021iterative} and residual-driven corrections of operators for stationary problems \citep{cao2023residual, jha2024residual}. None of them corrects a time-dependent operator under uncertainty.

The cost of time-stepping solvers lies in the integration, which is sequential over the horizon and repeated for every scenario and every iteration. A solution map removes it, since it returns the state at any time for all scenarios in one batched evaluation, and its gradient comes from automatic differentiation. We approximate this map by a physics-informed neural operator, the DeepONet \citep{lu2021deeponet, wang2021learning}, and train it on its residual alone, using \emph{generic controls} and no trajectory data, so a new target, horizon, or cost does not require retraining, as long as the dynamics and scenarios stay within the trained range. However, the operator accumulates error over the horizon, leading to large errors in the end position, which motivates our correction to its prediction.

The error between the prediction and the true state obeys, to first order, a linear differential equation driven by the residual \citep{hillebrecht2022certified}. By propagating the residual through the linearized dynamics, the first-order error can be removed, a technique we refer to as the \emph{correction}. The correction is, effectively, a single Newton step on the dynamics \citep{bellman1965quasilinearization}, started from the operator as in \citet{cao2023residual} and evaluated in parallel as in \citet{lim2024parallelizing}, and it follows the defect correction of ODE numerics \citep{zadunaisky1976estimation, stetter1978defect, skeel1986thirteen}. Related are the stabilized parallel Newton iterations of \citet{gonzalez2024towards}, parallel scans in control \citep{sarkka2023temporal} and operator predictions as starts of Newton's method \citep{aghili2025accelerating}.  We carry the correction to time-dependent dynamics under uncertainty with explicit constants, a computable certificate, and a cost guarantee. Since every step of the correction is an affine map, all scenarios and time steps are corrected by one parallel prefix scan \citep{blelloch1990prefix_chapter, lim2024parallelizing}, making it efficient to compute inside the optimization loop.

We deliver two main contributions. First, we provide a residual-based error correction for a physics-informed neural operator, with its analysis and a computable certificate of the corrected error. Second, unlike related approaches, we apply the corrected operator to a long-horizon problem with uncertain dynamics, the control of a fixed-wing UAV, where the corrected operator reaches the root mean square (RMS) miss of the reference single-shooting solvers with Adam and Ipopt within 0.6 and 2.5 m, respectively, and solves the benchmark 1.2 to 8 times faster than the Adam solver and 15 to 29 times faster than Ipopt at equal sample size. Section~\ref{section:problem} states the problem, Section~\ref{sec:methods_and_analysis} presents the proposed method and analysis, Section~\ref{sec:exp} describes the experiments of the UAV control problem, and Section~\ref{sec:conclusion} concludes.

\section{Uncertain optimal control problem for a UAV}\label{section:problem}
We study general uncertain optimal control problems, in which the equations describing the dynamics include uncertain parameters. Let $\xi \in \Xi$ be a random vector with distribution $\mathbb{P}$ that collects the uncertain parameters and the uncertain start state (initial conditions) of the dynamics. We name a realization of $\xi$ a \emph{scenario}. The state $s(t) \in \mathbb{R}^{n_s}$ of scenario $\xi$ under the control $u \colon [0, t_f] \to \mathbb{R}^{n_u}$ solves the dynamics equation expressed as an ordinary differential equation:
\begin{equation}\label{eq:dynamics}
  \dot s(t) = f\big(s(t), u(t), \xi\big), \qquad s(0) = s_0(\xi), \qquad t \in [0, t_f].
\end{equation}
Here $f \colon \mathbb{R}^{n_s} \times \mathbb{R}^{n_u} \times \Xi \to \mathbb{R}^{n_s}$ is the right-hand side of the dynamics and a dot denotes the first time derivative, i.e., $\dot s = \mathrm{d}s/\mathrm{d}t$. We specify the uncertain optimal control problem as
\begin{equation}\label{eq:ocp}
\begin{aligned}
  \min_{u \in \mathcal{U}} \quad & J(u) := \rho\Big[\ell\big(s(t_f; u, \xi)\big)\Big] + c(u) \\
  \text{s.t.} \quad & \eqref{eq:dynamics} \ \text{for every } \xi \in \Xi, \qquad g\big(\mathbb{E}_\xi[\, s(t; u, \xi)\,]\big) \le 0 \ \text{ for } t \in [0, t_f],
\end{aligned}
\end{equation}
where $\ell$ is the terminal cost of one scenario, i.e., at the end time $t_f$, $c$ is a deterministic control cost, $g$ a path constraint on the mean state, and $\rho$ a risk measure that maps the random terminal cost to a number. We assume that the control is open-loop, meaning that one $u \in \mathcal{U}$ is applied to every scenario without state feedback. Further, the admissible set $\mathcal{U}$ collects the deterministic constraints.

In this paper $\rho$ is the expectation, $\rho[\,\cdot\,] = \mathbb{E}_\xi[\,\cdot\,]$, and a  coherent risk measure \citep{artzner1999coherent}, such as the conditional value-at-risk \citep{rockafellar2000optimization}, fits the same formulation. Problem~\eqref{eq:ocp} is solved as a \emph{sample-average approximation}, i.e., for $M$ independent scenarios $\xi_1, \dots, \xi_M \sim \mathbb{P}$,%
\begin{equation}\label{eq:saa}
  J_M(u) := \frac{1}{M} \sum_{i=1}^{M} \ell\big(s(t_f; u, \xi_i)\big) + c(u),
\end{equation}
and every evaluation of $J_M$ integrates \eqref{eq:dynamics} $M$ times. For a fixed control $u$ with $J(u) < \infty$, the strong law of large numbers gives $J_M(u) \to J(u)$ almost surely as $M \to \infty$ \citep{shapiro2009lectures}.

Our focus lies on the fixed-wing UAV problem of \citet[eqs.~31--41]{lambrianides2020scalable}. The state $s = (x, y, z, v, \gamma, \sigma)$ holds the position $p = (x, y, z)$, the speed $v$, the flight-path angle $\gamma$ and the heading $\sigma$. The inputs are the thrust $T$, the angle of attack $\alpha$ and the bank angle $\mu$. The point-mass flight equations are
\begin{equation}\label{eq:uav}
\begin{aligned}
  \dot p &= v\,\big(\cos\gamma\cos\sigma,\ \cos\gamma\sin\sigma,\ \sin\gamma\big)^{\!\top}, &
  \dot v &= \big(T\cos\alpha - F_D\big)/m_a - g_0\sin\gamma, \\
  \dot\gamma &= \big(F_L + T\sin\alpha\big)\cos\mu/(m_a v) - g_0\cos\gamma/v, &
  \dot\sigma &= \big(F_L + T\sin\alpha\big)\sin\mu/(m_a v\cos\gamma),
\end{aligned}
\end{equation}
with mass $m_a$, gravity $g_0$, lift and drag $F_L = \tfrac12 \varrho(z)\, v^2 S\, C_L$ and $F_D = \tfrac12 \varrho(z)\, v^2 S\, C_D$, wing area $S$, an exponential air-density profile $\varrho(z) = 1.21\, e^{-z/8000}\ \text{kg/m}^3$, and the aerodynamic coefficients
\begin{equation}\label{eq:aero}
  C_L = c_a \sin\alpha - c_z \cos\alpha, \qquad C_D = -c_a \cos\alpha - c_z \sin\alpha,
\end{equation}
with $c_a = C_{x0} + C_{x\alpha}\,\alpha$ and $c_z = C_{z0} + C_{z\alpha}\,\alpha$. The control is the input vector $u = (T, \alpha, \mu)$, which is the same for every scenario, starts from a given value $u_0$ and has bounded rates. The aircraft has to reach a target point $p^*$ at the final time $t_f$. In the notation of \eqref{eq:ocp},
\begin{equation}\label{eq:uav_cost}
  \ell\big(s(t_f)\big) = \big\| p(t_f) - p^* \big\|^2, \qquad
  c(u) = \lambda_u \int_0^{t_f} \|\dot u(t)\|^2 \, dt, \qquad
  \rho = \mathbb{E}_\xi ,
\end{equation}
where $\|\cdot\|$ is the Euclidean norm. The state constraint $g$ of \eqref{eq:ocp} bounds the mean speed and the mean flight-path angle,
\begin{equation}\label{eq:uav_bounds}
  v_{\min} \le \frac{1}{M}\sum_{i=1}^{M} v(t; u, \xi_i) \le v_{\max}, \qquad
  \Big| \frac{1}{M}\sum_{i=1}^{M} \gamma(t; u, \xi_i) \Big| \le \gamma_{\max}, \qquad t \in [0, t_f].
\end{equation}
A scenario $\xi$ has ten components, six uncertain start states $s_0(\xi) = (p_0, v_0, \gamma_0, \sigma_0)$ and four uncertain aerodynamic coefficients $C_{x0}$, $C_{x\alpha}$, $C_{z0}$ and $C_{z\alpha}$ of \eqref{eq:aero}, which set the lift and the drag. All values used in the experiments are listed in Table~\ref{tab:setup} in Appendix~\ref{app:problem}.

\section{Method and analysis}\label{sec:methods_and_analysis}

\subsection{Neural operator}
\label{sec:operator}

For solving the uncertain optimal control problem, we replace the dynamics equation \eqref{eq:dynamics} with its solution map $s \colon [0, t_f]\times\mathcal{U}\times\Xi\to\mathbb{R}^{n_s}$, which we approximate by a neural operator $\hat s$. For a control $u$, a scenario $\xi$, and any time $t \in [0, t_f]$, it returns a prediction $\hat{s}(t; u, \xi) \in \mathbb{R}^{n_s}$ of the state in one network evaluation. We use the branch--trunk form of a DeepONet \citep{lu2021deeponet}. A branch network $b$ sees a discretization of the control on $[0, t_f]$, and the trunk $\psi$ sees the time $t$ and the scenario $\xi$. Formally, this can be written as
\begin{equation}\label{eq:operator}
\hat s_k(t; u, \xi) = s_{0,k}(\xi) + \frac{t}{t_f}\, d_k \Big( \sum_{j=1}^{n_\psi} b_{kj}(u)\, \psi_{kj}(t, \xi) + \beta_k \Big), \qquad k = 1, \dots, n_s,
\end{equation}
with $n_\psi$ terms per state component, state scales $d_k$ and biases $\beta_k$. The factor $t/t_f$ makes $\hat s(0; u, \xi) = s_0(\xi)$ hold exactly. Further, this form has two beneficial consequences for the optimization solver. First, the trunk does not see the control, so its value and exact time derivatives are computed once, before the optimization. Second, the per-iteration cost of the operator is then one branch evaluation and one tensor contraction of the stored trunk values.

For physics-informed neural operators, the \emph{residual}
\begin{equation}\label{eq:residual}
  R(t; u, \xi) := \partial_t \hat s(t; u, \xi) - f\big(\hat s(t; u, \xi), u(t), \xi\big)
\end{equation}
measures how much the prediction violates the dynamics and is computable from the operator and $f$ alone. The operator is trained by minimizing the residual over its weights \citep{wang2021learning, li2024physics, goswami2023physics}.

\subsection{Single shooting solution}\label{sec:single_shooting}

By discretizing the control and replacing the dynamics with the neural operator, the sample-average problem \eqref{eq:saa}--\eqref{eq:uav_bounds} becomes a \emph{single-shooting} problem as the states are eliminated by the forward map, and the only unknown is the control \citep{biegler2010nonlinear}. To obtain smooth controls, we parametrize $u$ by cubic B-splines with $n_b$ coefficients per control channel,
\begin{equation}\label{eq:spline}
  u(t; w) = u_{\mathrm{mid}} + u_{\mathrm{half}} \odot \sum_{k=1}^{n_b} w_k\, B_k(t),
\end{equation}
where $\odot$ is the component-wise product, $u_{\mathrm{mid}}$ and $u_{\mathrm{half}}$ are the centre and the half-width of the admissible box defined by $\mathcal{U}$, and $w \in \mathcal{W} := [-1, 1]^{n_u \times n_b}$ are the weights. The B-splines are nonnegative and sum to one, so $u(t; w)$ stays in the box for every $t$ and every $w \in \mathcal{W}$. The raw surrogate cost is then
\begin{equation}\label{eq:rawcost}
  \hat J_M(u) := \frac{1}{M} \sum_{i=1}^{M} \ell\big(\hat s(t_f; u, \xi_i)\big) + c(u), \qquad u = u(\cdot\,; w).
\end{equation}
The raw cost \eqref{eq:rawcost} can be minimized directly by a gradient-based optimizer, with gradients obtained using automatic differentiation through the network branch and the sum \eqref{eq:operator}. However, the optimizer then sees only the predicted end state $\hat s(t_f; u, \xi)$ and is free to exploit any error it can \citep{lundqvist2025employing}.
Therefore, we correct the predicted state using an estimate of its error, so that the cost and its gradient follow a state that is sufficiently accurate.

\subsection{Error correction}\label{sec:correction}
We define the \emph{error} between the true solution map $s$ and our neural operator map $\hat s$ as $e(t; u, \xi) := \hat s(t; u, \xi) - s(t; u, \xi)$. Hereinafter, we drop the arguments $u$ and $\xi$ when they are clear. Further, we write $\|\cdot\|$ for a norm on $\mathbb{R}^{n_s}$ with its induced operator norms, the Euclidean norm unless stated otherwise. The error of the operator obeys a differential equation of its own. Differentiating $e = \hat s - s$ with respect to $t$ and expanding $f(s)$ about the prediction $\hat{s}$ by Taylor's theorem gives
\begin{equation}\label{eq:error_ode}
  \dot e(t) = \hat A(t)\, e(t) + R(t) - q(t), \qquad e(0) = 0, \qquad \hat A(t) := \partial_s f\big(\hat s(t), u(t), \xi\big),
\end{equation}
where $\hat A \in \mathbb{R}^{n_s \times n_s}$ is the Jacobian of the dynamics along the prediction $\hat{s}$ and $q(t) = \int_0^1 (1-\theta)\, \Big( e(t)^\top\, \partial_s^2 f_i\big(\hat s(t) - \theta\, e(t), u(t), \xi\big)\, e(t) \Big)_{i=1}^{n_s}\, d\theta$ is the Taylor remainder, which is of second order in $e$. The residual acts as a forcing on the error, which the linearized dynamics $\hat A$ propagate over the horizon. Thus, for the UAV, a small residual in the angles accumulates into a large error in the end position. Everything in \eqref{eq:error_ode} except $q$ is computable from the operator and $f$. Dropping $q$ gives the linear equation
\begin{equation}\label{eq:estimate_ode}
  \dot{\hat e}(t) = \hat A(t)\, \hat e(t) + R(t), \qquad \hat e(0) = 0,
\end{equation}
whose solution is the residual propagated by the linearized dynamics,
\begin{equation}\label{eq:estimate}
  \hat e(t) = \int_0^t \hat\Phi(t, \tau)\, R(\tau)\, d\tau, \qquad
  \partial_t \hat\Phi(t, \tau) = \hat A(t)\, \hat\Phi(t, \tau), \quad \hat\Phi(\tau, \tau) = I,
\end{equation}
where $\hat\Phi(t, \tau) \in \mathbb{R}^{n_s \times n_s}$ is the state transition matrix of \eqref{eq:estimate_ode} \citep{brockett2015finite}. The estimate $\hat e$ agrees with $e$ to first order, since the same formula applied to \eqref{eq:error_ode} gives $e(t) = \int_0^t \hat\Phi(t, \tau)\,[R(\tau) - q(\tau)]\, d\tau$ and the two differ by the propagated remainder.

As in defect correction \citep{stetter1978defect}, we define the \emph{corrected state} $\bar s$ as
\begin{equation}\label{eq:corrected}
  \bar s(t; u, \xi) := \hat s(t; u, \xi) - \hat e(t; u, \xi),
\end{equation}
whose error $\bar s - s$ is of second order in $e$. Alternatively derived, the correction is one Newton step for the root of $s \mapsto \dot s - f(s)$, started at $\hat s$ \citep{bellman1965quasilinearization, lim2024parallelizing}. The residual plus the linearization at $\hat s$, applied to the step, vanishes, $R + \tfrac{d}{dt}(\bar s - \hat s) - \hat A\,(\bar s - \hat s) = 0$ with $(\bar s - \hat s)(0) = 0$, which gives the same differential equation \eqref{eq:estimate_ode} by writing $\bar s - \hat s = -\hat e$. The following lemma states the quadratic convergence of this Newton step.

\begin{lemma}[second-order correction]\label{lem:identity}
Let $f$ be twice continuously differentiable in $s$ on convex sets $\mathcal C(t)$ that contain $s(t)$ and $\hat s(t)$ for all $t \in [0, t_f]$, let $\hat s(0) = s_0(\xi)$, and let $\hat K(t)$ be the propagated curvature \eqref{eq:K} on these sets. Then for every $t \in [0, t_f]$
\begin{align}
  \big\| \bar s(t) - s(t) \big\| &\le \tfrac12\, \hat K(t) \sup_{\tau \le t} \|e(\tau)\|^2, \label{eq:correction} \\
  \|e(t)\| &\le \|\hat e(t)\| + \tfrac12\, \hat K(t) \sup_{\tau \le t} \|e(\tau)\|^2 . \label{eq:error_bound}
\end{align}
\end{lemma}

The constant $\hat K(t)$ is the largest curvature of the dynamics between prediction and truth, propagated by the linearized dynamics and integrated up to time $t$, with the formula \eqref{eq:K} in Appendix~\ref{app:proofs}. Notice that Lemma~\ref{lem:identity} is an a priori statement. Its right-hand side contains the true error $e$, which is unknown, so it only gives the order of the correction. We obtain a computable bound if we use the corrected state $\bar s$ and bound its distance to the true state $s$ by its own residual $\bar R := \dot{\bar s} - f(\bar s) = f(\hat s) - \hat A\, \hat e - f(\bar s)$, which follows from \eqref{eq:residual} and \eqref{eq:estimate_ode}. Thus, let $\bar A := \partial_s f(\bar s, u, \xi)$ be the Jacobian at the corrected state and $\bar\Phi$ the transition matrix as in \eqref{eq:estimate}, as $\hat\Phi$ is for $\hat A$. We measure the residual $\bar R$ propagated along $\bar s$ by
\begin{equation}\label{eq:abar}
  \bar a := \sup_{t \le t_f} \int_0^t \big\| \bar\Phi(t, \tau)\, \bar R(\tau) \big\|\, d\tau .
\end{equation}
Finally, let $\bar K$ be the propagated curvature \eqref{eq:K_bar} along $\bar s$ on the tubes of radius $2\bar a$ around it, taken at its largest over $t$. The following theorem turns these two computable quantities into a bound on the error of the corrected state.

\begin{theorem}[certificate at the corrected state]\label{thm:certificate}
Let $f$ be twice continuously differentiable in $s$ on the balls $\bar{\mathcal C}(\tau) := \{\zeta \in \mathbb{R}^{n_s} \colon \|\zeta - \bar s(\tau)\| \le 2\bar a\}$, $\tau \in [0, t_f]$, which form a tube of radius $2\bar a$ around $\bar s$. If
\begin{equation}\label{eq:certificate}
  2\, \bar K\, \bar a < 1,
\end{equation}
then the true state exists on $[0, t_f]$ and $\|\bar s(t) - s(t)\| \le 2\bar a$ for all $t \in [0, t_f]$.
\end{theorem}
Theorem~\ref{thm:certificate} bounds the distance of the corrected state from the true state by $2\bar a$ whenever the certificate \eqref{eq:certificate} holds, and both $\bar a$ and $\bar K$ are computed from the operator and $f$ at the returned control, with no true trajectory. Condition \eqref{eq:certificate} is a Newton--Kantorovich type condition in affine covariant form \citep{kantorovich1982functional, deuflhard2004newton}. Evaluated at the prediction, i.e., with $\hat s$, $R$, $\hat\Phi$ and $\hat K$ in place of $\bar s$, $\bar R$, $\bar\Phi$ and $\bar K$, condition \eqref{eq:certificate} failed for every scenario in our experiments, emphasizing the need for the correction. We therefore evaluate the certificate at the corrected state. The proofs and the formulas are in Appendix~\ref{app:proofs}.

Although every quantity can be computed from the operator, the state transition matrix $\hat\Phi$ has no closed form for a time-varying Jacobian. For a constant Jacobian $\hat A$ the transition matrix equals $\exp\big((t - \tau)\hat A\big)$. Thus, on a small enough grid spacing, the transition matrix can be approximated by that of a constant Jacobian. Therefore, we discretize \eqref{eq:estimate_ode} on a grid $t_k = kh$, $k = 0, \dots, N$, by the one-step scheme
\begin{equation*}
  \hat e(t_{k+1}) = \hat\Phi_k\, \hat e(t_k) + \hat r_k, \qquad
  \hat\Phi_k = I + h\hat A_{k+1/2} + \tfrac12 (h\hat A_{k+1/2})^2, \qquad
  \hat r_k = \tfrac{h}{2}\big(\hat\Phi_k R(t_k) + R(t_{k+1})\big),
\end{equation*}
with $\hat A_{k+1/2} = \tfrac12\big(\hat A(t_k) + \hat A(t_{k+1})\big)$, which is second order in $h$. Each step is an affine map, and affine maps compose associatively, i.e., $(\hat\Phi_b, \hat r_b)\circ(\hat\Phi_a, \hat r_a) = (\hat\Phi_b\hat\Phi_a,\ \hat\Phi_b \hat r_a + \hat r_b)$. Therefore, all $\hat e(t_k)$ of all scenarios follow from one parallel prefix scan in $O(\log N)$ stages \citep{blelloch1990prefix_chapter,lim2024parallelizing} instead of $N$ sequential steps, which keeps the computation at a fraction of the cost of a time-stepping integration. The residual, the Jacobian, and the scan are differentiated (with automatic differentiation) together with the rest of the cost, and as such, the optimizer only sees the corrected states. With the correction, the solver minimizes the corrected cost
\begin{equation}\label{eq:checkedcost}
  \bar J_M(u) := \frac{1}{M} \sum_{i=1}^{M} \ell\big(\bar s(t_f; u, \xi_i)\big) + c(u), \qquad u = u(\cdot\,; w)
\end{equation}
in place of \eqref{eq:rawcost}. The correction is applied in every iteration for the fixed set of $M$ sampled scenarios $\xi_1, \dots, \xi_M$ chosen once.

\subsection{Guarantees for the corrected problem}
The solver minimizes the corrected cost $\bar J_M$ in place of the true cost $J$. The gap $|J(u) - \bar J_M(u)|$ has two sources. The first is the sampling error $|J(u) - J_M(u)|$ of the $M$ scenarios. The second is the mean error $\bar\epsilon(u) := \frac{1}{M} \sum_{i=1}^{M} \| \bar s_i(t_f; u) - s_i(t_f; u) \|$ of the corrected end states, which is of second order in the error of the prediction by Lemma~\ref{lem:identity}. The following theorem bounds the distance between the two costs, $J$ and $\bar J_M$, by these two terms for every admissible control, so the optimizer can only be misled where the surrogate cost deviates from the true one.

\begin{theorem}[error and suboptimality of the corrected cost]\label{thm:honest}
Let the assumptions of Lemma~\ref{lem:identity} hold for every sampled scenario, and let $\ell$ be Lipschitz continuous with constant $L_\ell$ on a set that contains the true and the corrected end states of every scenario for all $u \in \mathcal{U}$. Then, for every $u \in \mathcal{U}$
\begin{equation}\label{eq:honest}
  \big| J(u) - \bar J_M(u) \big| \le \big| J(u) - J_M(u) \big| + L_\ell\, \bar\epsilon(u),
\end{equation}
and $\big\| \bar s_i(t_f; u) - s_i(t_f; u) \big\| \le 2 \bar a_i$ for every scenario $i$ that satisfies \eqref{eq:certificate} at $u$. If $\hat u$ minimizes $\bar J_M$ over $\mathcal{U}_{\mathcal{W}} := \{u(\cdot\,; w) \colon w \in \mathcal{W}\}$ up to a tolerance $\varepsilon_{\mathrm{opt}}$ and $u^*$ minimizes $J$ over $\mathcal{U}_{\mathcal{W}}$, then
\begin{equation}\label{eq:subopt}
  J(\hat u) - J(u^*) \le \varepsilon_{\mathrm{opt}} + 2 \sup_{u \in \mathcal{U}_{\mathcal{W}}} \big| J(u) - J_M(u) \big| + L_\ell \big( \bar\epsilon(\hat u) + \bar\epsilon(u^*) \big).
\end{equation}
\end{theorem}

Theorem~\ref{thm:honest} says that the solver can be misled only by the sampling error and by the error of the corrected end states, and it says so for every control, not only at the optimal solution. In \eqref{eq:subopt}, the sampling error enters as a supremum over all controls because $\hat u$ is selected with the same $M$ scenarios that define $J_M$, so the law of large numbers for a fixed control does not apply to it. If $\Xi$ is compact and the true state exists for every $u \in \mathcal U_{\mathcal W}$ and $\xi \in \Xi$, the uniform law of large numbers \citep[Theorem~7.48]{shapiro2009lectures} makes the supremum vanish with probability one as $M$ grows. Theorem~\ref{thm:certificate} bounds the correction error at each control for every certified scenario, and it vanishes with the residual. A gradient method returns a local minimizer, for which $\varepsilon_{\mathrm{opt}}$ can be large, and \eqref{eq:honest} still holds. The proof is in Appendix~\ref{app:proofs}.

\section{Experiments}\label{sec:exp}
\subsection{Problem, setup and reference solutions}

We apply the corrected operator to the UAV problem of Section~\ref{section:problem}. The corrected forward map $\bar s(\cdot\,; u, \xi) = \hat s - \hat e$ of \eqref{eq:corrected} replaces the dynamics stated in \eqref{eq:dynamics}, so with the spline inputs \eqref{eq:spline} and $\ell$ and $c$ from \eqref{eq:uav_cost}, the operator solver solves
\begin{equation}\label{eq:solver_problem}
  \min_{w \in \mathcal W}\ \bar J_M\big(u(\cdot\,; w)\big) \qquad \text{s.t.} \quad s_i = \bar s(\cdot\,; u(\cdot\,; w), \xi_i),\ i = 1, \dots, M, \quad \text{and } \eqref{eq:uav_bounds}.
\end{equation}
Appendix~\ref{app:problem} lists the problem data in Table~\ref{tab:setup}, which are those of \citet{lambrianides2020scalable}.

Each input uses a cubic B-spline with $n_b = 26$ coefficients and knots every 2.2 s. The first coefficient of each input is fixed so that $u(0) = u_0$, which leaves 75 free coefficients, the same for every solver. The operator of Appendix~\ref{app:details} is trained once on the residual \eqref{eq:residual} of 60{,}000 generic controls, without trajectory data, in 2.4 h on one GPU, and the correction runs on all $M$ scenarios on the grid $h = 1$ s with $N = 50$ steps.

We compare the corrected neural operator against two \emph{reference solvers}, which solve \eqref{eq:solver_problem} with the states of a fourth-order Runge--Kutta (RK4) integration at step 0.5 s in place of the corrected map, on the same splines and the same $M$ scenarios. The first uses the same Adam schedule, batched over the scenarios, compiled with XLA in JAX and runs on the GPU. A larger step is not possible, as it violates the mean-state bounds (see Appendix~\ref{app:rk4steps}). The second solves the problem using the same RK4 states as a nonlinear program with hard constraints, implemented in CasADi \citep{andersson2019casadi} and Ipopt \citep{wachter2006implementation} with exact Hessians, which are 11 to 12 times faster than a limited-memory Hessian here. We run the Ipopt solver on the CPU as CasADi has no GPU backend. We estimate that Ipopt with GPU evaluations would reach the time of the reference Adam solver at large $M$. Therefore, we use the Ipopt solver mainly as a reference for solution quality, since its solve times are on different hardware. The operator solver and Adam RK4 minimize the cost with logarithmic barriers and quadratic penalties on the input, rate and mean-state bounds, by Adam with 20{,}000 steps of size $3 \cdot 10^{-3}$, cosine decay, and projection onto $\mathcal W$. Thus the operator solver and the reference Adam solver are identical except for how they solve the dynamics.

Appendix~\ref{app:problem} states both problems in detail, and Appendix~\ref{app:details} provides further details on the solvers. The corrected neural operator and the reference Adam solver run on one GPU (RTX 5070) and Ipopt on the CPU (10 cores of an Intel Core Ultra 5 225). The solve time is the time of the optimization loop, without the compilation (7 s for the operator, 40 to 53 s for the reference Adam solver), the precomputation of the trunk and the build of the nonlinear program (10 to 102 s).

\subsection{Results}

\begin{figure}[t]\centering
\includegraphics[width=0.9\textwidth]{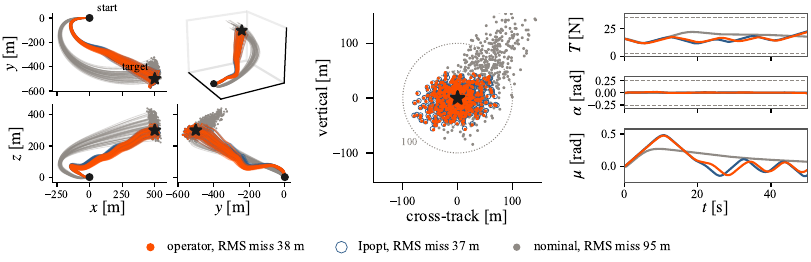}
\caption{Left, 60 of 500 held-out flights. Middle, all 500 endpoints seen along the mean approach direction. Right, the inputs with their bounds. The legend gives the RMS miss.}\label{fig:traj}
\end{figure}

For every solution, we report the \emph{RMS miss} $\big(\frac{1}{N_{\mathrm{test}}}\sum_{i=1}^{N_{\mathrm{test}}} \|p_i(t_f) - p^*\|^2\big)^{1/2}$ over $N_{\mathrm{test}}=200{,}000$ held-out scenarios flown with RK4 at step 0.1 s, which we refer to as the \emph{realized RMS miss}. The \emph{apparent RMS miss} is the RMS miss a solver computes on its own $M$ scenarios, and the \emph{model gap} is its difference from the RMS miss of the same control flown with RK4 at 0.1 s on those $M$ scenarios. For every problem, we draw three sets of scenarios with different random seeds and report the solve times as means across the seeds. Figure~\ref{fig:traj} shows a solution for the second target $p^* = (500, -500, 300)$ m. The operator solution and the Ipopt solution have almost the same controls, with thrusts that differ by at most 2.8 N, and RMS misses of 38 and 37 m on 500 held-out scenarios. At the operator solution, the correction reduces the RMS endpoint error of the operator from 76 m to 4.1 m. The solution of the nominal problem, which is optimal for the mean scenario only, has an RMS miss of 95 m, since the mean of its endpoints lies 67 m from the target.

\begin{table}[t]\centering\footnotesize
\caption{Benchmark at equal $M$. RMS miss over 200{,}000 held-out scenarios, largest model gap over the seeds, mean solve time, and the factor of the solve times against the operator.}\label{tab:results}
\setlength{\tabcolsep}{2pt}
\begin{tabular*}{\textwidth}{@{\extracolsep{\fill}}rccccccccc@{}}
\toprule
 & \multicolumn{3}{c}{operator with the correction} & \multicolumn{2}{c}{reference RK4, Adam} & \multicolumn{2}{c}{reference RK4, Ipopt} & \multicolumn{2}{c}{factor} \\
\cmidrule(lr){2-4}\cmidrule(lr){5-6}\cmidrule(lr){7-8}\cmidrule(lr){9-10}
$M$ & RMS miss [m] & model gap [m] & solve [s] & RMS miss [m] & solve [s] & RMS miss [m] & solve [s] & Adam & Ipopt \\
\midrule
128 & 52.3/52.4/52.9 & 0.1 & 9.7 & 52.5/52.7/53.0 & 81 & 51.8/51.7/52.5 & 142 & 8$\times$ & 15$\times$ \\
400 & 52.5/52.4/52.3 & $<0.1$ & 21.2 & 52.6/52.6/52.4 & 161 & 51.7/51.8/51.6 & 413 & 8$\times$ & 20$\times$ \\
1600 & 52.0/52.0/52.2 & $<0.1$ & 58.5 & 52.3/52.1/52.2 & 200 & 51.7/51.6/51.7 & 1635 & 3$\times$ & 28$\times$ \\
4800 & 52.1/52.0/52.0 & $<0.1$ & 161 & 52.2/52.3/52.4 & 249 & 51.6/51.6/51.6 & 4624 & 1.5$\times$ & 29$\times$ \\
9600 & 52.0/52.0/52.0 & $<0.1$ & 334 & 52.2/52.1/52.3 & 412 & -- & -- & 1.2$\times$ & -- \\
\bottomrule
\end{tabular*}

\end{table}

\begin{figure}[t]\centering
\includegraphics[width=0.9\textwidth]{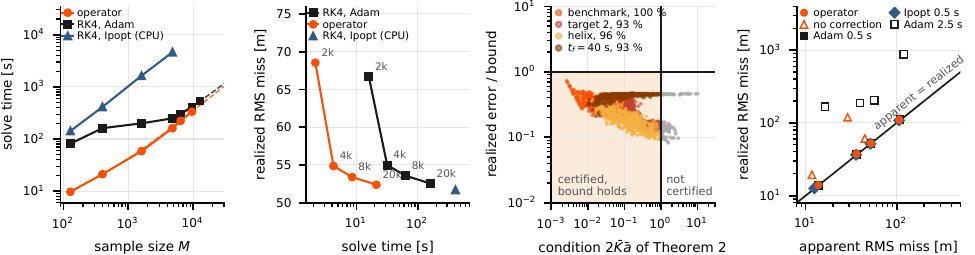}
\caption{Left to right. Solve time against $M$, dashed the linear cost. RMS miss against solve time at $M = 400$. Certificate of Theorem~\ref{thm:certificate}, seed 0. Apparent against realized miss.}\label{fig:results}
\end{figure}

Table~\ref{tab:results} and Figure~\ref{fig:results} compare the solvers at equal sample size. For every $M$ and every seed, the operator reaches the realized RMS miss of the reference Adam solver within 0.4 m, and Ipopt reaches a 0.4 to 1.0 m lower miss than both Adam solvers, since it enforces the mean-state bound as a hard constraint with an exact Hessian, while the barrier and the fixed budget of the Adam loop leave this margin. The operator solves the problem in 10 s at $M = 128$ and in 21 s at $M = 400$, which is 8 times faster than the reference Adam solver and 15 and 20 times faster than Ipopt. Above a few thousand scenarios, both costs grow at about 1.8 $\mu$s per scenario and iteration, so the factor against the Adam solver falls slowly: 3 at $M = 1600$, 1.5 at 4800 and 1.2 at 9600, the largest sample size at which the operator fits the 12 GB of the GPU (Figure~\ref{fig:results}, left). Both Adam solvers converge in the same number of steps, and each step of the operator is 8 times cheaper, so at $M = 400$ the operator reaches every level of the realized RMS miss 8 times sooner (Figure~\ref{fig:results}, second panel). The third panel of Figure~\ref{fig:results} evaluates the certificate at the returned controls of the four tasks for seed 0, in the norm scaled by the state scales $d_k$ of \eqref{eq:operator}, with $\bar s$ on a grid of 0.025 s. The condition \eqref{eq:certificate} holds for 88 to 100\,\% of 500 held-out scenarios, and the realized error of the corrected state stays below the bound $2\bar a$ in every certified scenario. Appendix~\ref{app:further} gives the cost of each part of one iteration and the realized RMS miss after shorter Adam budgets and with a second correction step.

\begin{table}[t]\centering\footnotesize
\caption{Four tasks at $M = 400$, three seeds, misses in metres as in Table~\ref{tab:results}, $\|\bar p - p\|$ on 500 held-out scenarios, largest in brackets. A dagger marks a mean flight-path angle above $\gamma_{\max}$.}\label{tab:tasks}
\setlength{\tabcolsep}{1.5pt}\scriptsize\renewcommand{\arraystretch}{0.92}
\begin{tabular*}{\textwidth}{@{\extracolsep{\fill}}lccccccc@{}}
\toprule
 & \multicolumn{3}{c}{operator with the correction} & Adam & Ipopt & no correction & RK4 2.5 s \\
\cmidrule(lr){2-4}\cmidrule(lr){5-5}\cmidrule(lr){6-6}\cmidrule(lr){7-7}\cmidrule(lr){8-8}
task & RMS miss & model gap & $\|\bar p - p\|$ & RMS miss & RMS miss & apparent/realized & apparent/realized \\
\midrule
benchmark & 52.5/52.4/52.3 & $<0.1$ & 2.1 (4.7) & 52.6/52.6/52.4 & 51.7/51.8/51.6 & 45 / 61$^{\smash{\dagger}}$ & 57 / 206$^{\smash{\dagger}}$ \\
target 2 & 37.6/37.5/37.8 & 0.5 & 4.1--5.9 (11.5) & 37.2/37.1/37.2 & 36.6/36.6/36.6 & 29 / 119$^{\smash{\dagger}}$ & 40 / 188$^{\smash{\dagger}}$ \\
helix & 14.1/14.0/14.0 & 0.1 & 1.4--1.5 (4.9) & 14.0/14.0/13.9 & 12.7/12.7/12.6 & 12 / 19$^{\smash{\dagger}}$ & 16 / 167$^{\smash{\dagger}}$ \\
$t_f = 40$ s & 110.6/110.8/110.2 & 0.8 & 2.5--2.6 (26) & 110.8/110.6/110.5 & 108.6/108.3/108.3 & -- & 120 / 878$^{\smash{\dagger}}$ \\
\bottomrule
\end{tabular*}
\end{table}

Table~\ref{tab:tasks} compares five solvers on four tasks with the same operator, trained once. The tasks are the benchmark $p^* = (500, 500, 500)$ m, a second target $p^* = (500, -500, 300)$ m (target 2), a third target $p^* = (0, 0, 400)$ m above the start, which the aircraft reaches on a climbing helix (helix), and the shorter horizon $t_f = 40$ s for the benchmark, for which the inputs are held constant after 40 s on the branch grid. The operator with the correction reaches the realized RMS miss of the reference Adam solver within 0.6 m on every task and trails Ipopt by 0.5 to 2.5 m, as the Adam solver does. The operator solves each task in 21 to 29 s, the reference Adam solver in 131 to 161 s, and Ipopt in 161 to 628 s. Without the correction, on three tasks, and with the reference solver at step 2.5 s, on all four, the optimizer exploits the model error, and the apparent RMS misses lie far below the realized ones (Figure~\ref{fig:results}, right), with mean flight-path angles above the bound.

\section{Conclusion}\label{sec:conclusion}

We showed that a learned solution map can replace the time-stepping integrator inside the optimizer if it is corrected by removing the first-order error. The correction is one Newton step on the dynamics, computed with a parallel prefix scan. The scan costs $O(n_s^3)$ per step, so systems with many states need banded or low-rank Jacobians. Each iteration also computes the residual and its Jacobian on the correction grid. For the corrected state, we provided a computable certificate that bounds the error between the corrected state and the true state. We further showed that the gap between the true cost and the corrected sample-average cost is bounded by the sampling error and the same state error. We compared the operator with a batched RK4 Adam solver on the same GPU and with Ipopt on the CPU as the reference for accuracy. The operator is 8 times faster than the batched integrator for a few hundred scenarios on one GPU, and 1.2 times faster at 9600, where it reaches the GPU's memory limit, so larger sample sizes require compressed or streamed features.  A larger RK4 step would make the integrator faster, but at step sizes of 1.0 and 2.5 s its solutions violate the bound on the mean flight-path angle (Appendix~\ref{app:rk4steps}). Further, the correction can start from the previous iterate of the optimizer instead of the operator prediction, which removes the operator from the optimization loop. In preliminary runs, this variant matched the operator’s time on the benchmark but required a step cap and periodic exact RK4 rollouts to converge, and its realized RMS miss was still 7 to 32 m larger on the 40 s task, so we chose to leave it out of the study.

Finally, our work has limitations. Our results bound the cost but do not provide guarantees on the convergence of the gradient method to a stationary point of the true cost, and the B-spline parametrization excludes bang-bang controls. Further assessments that consider coherent risk measures, systems with more states and time-varying disturbances such as wind gusts are possible directions for future research.

\subsection*{Generative AI Use}
Generative AI tools (Claude, Anthropic) assisted with the literature search, derivations, coding, experiments and editing. The authors verified all AI-assisted content and take responsibility for it.

\acks{This work was supported by the Finnish Ministry of Education and Culture’s Pilot for Doctoral Programmes (Pilot project Mathematics of Sensing, Imaging and Modelling).}

\bibliography{main}

\begin{thebibliography}{49}
\providecommand{\natexlab}[1]{#1}
\providecommand{\url}[1]{\texttt{#1}}
\expandafter\ifx\csname urlstyle\endcsname\relax
  \providecommand{\doi}[1]{doi: #1}\else
  \providecommand{\doi}{doi: \begingroup \urlstyle{rm}\Url}\fi

\bibitem[Aghili et~al.(2025)Aghili, Franck, Hild, Michel-Dansac, and Vigon]{aghili2025accelerating}
Joubine Aghili, Emmanuel Franck, Romain Hild, Victor Michel-Dansac, and Vincent Vigon.
\newblock Accelerating the convergence of {Newton}'s method for nonlinear elliptic {PDEs} using {Fourier} neural operators.
\newblock \emph{Communications in Nonlinear Science and Numerical Simulation}, 140:\penalty0 108434, 2025.

\bibitem[Allen et~al.(2022)Allen, Lopez-Guevara, Stachenfeld, Sanchez-Gonzalez, Battaglia, Hamrick, and Pfaff]{allen2022physical}
Kelsey~R. Allen, Tatiana Lopez-Guevara, Kimberly Stachenfeld, Alvaro Sanchez-Gonzalez, Peter Battaglia, Jessica Hamrick, and Tobias Pfaff.
\newblock {Physical design using differentiable learned simulators}.
\newblock arXiv:2202.00728, 2022.

\bibitem[Andersson et~al.(2019)Andersson, Gillis, Horn, Rawlings, and Diehl]{andersson2019casadi}
Joel A.~E. Andersson, Joris Gillis, Greg Horn, James~B. Rawlings, and Moritz Diehl.
\newblock {CasADi: A software framework for nonlinear optimization and optimal control}.
\newblock \emph{Mathematical Programming Computation}, 11\penalty0 (1):\penalty0 1--36, 2019.
\newblock \doi{10.1007/s12532-018-0139-4}.

\bibitem[Artzner et~al.(1999)Artzner, Delbaen, Eber, and Heath]{artzner1999coherent}
Philippe Artzner, Freddy Delbaen, Jean-Marc Eber, and David Heath.
\newblock {Coherent Measures of Risk}.
\newblock \emph{Mathematical Finance}, 9\penalty0 (3):\penalty0 203--228, 1999.
\newblock \doi{10.1111/1467-9965.00068}.

\bibitem[Bellman and Kalaba(1965)]{bellman1965quasilinearization}
Richard~E. Bellman and Robert~E. Kalaba.
\newblock \emph{Quasilinearization and Nonlinear Boundary-Value Problems}.
\newblock American Elsevier, New York, 1965.

\bibitem[Biegler(2010)]{biegler2010nonlinear}
Lorenz~T. Biegler.
\newblock \emph{{Nonlinear Programming: Concepts, Algorithms, and Applications to Chemical Processes}}.
\newblock MOS-SIAM Series on Optimization. Society for Industrial and Applied Mathematics, Philadelphia, 2010.
\newblock \doi{10.1137/1.9780898719383}.

\bibitem[Blelloch(1990)]{blelloch1990prefix_chapter}
Guy~E. Blelloch.
\newblock Prefix sums and their applications.
\newblock In John~H. Reif, editor, \emph{Synthesis of Parallel Algorithms}. Morgan Kaufmann Publishers, 1990.

\bibitem[Brockett(2015)]{brockett2015finite}
Roger~W. Brockett.
\newblock \emph{{Finite Dimensional Linear Systems}}.
\newblock Classics in Applied Mathematics. Society for Industrial and Applied Mathematics, Philadelphia, 2015.
\newblock \doi{10.1137/1.9781611973884}.
\newblock Reprint of the 1970 Wiley edition.

\bibitem[Cao et~al.(2023)Cao, O'Leary-Roseberry, Jha, Oden, and Ghattas]{cao2023residual}
Lianghao Cao, Thomas O'Leary-Roseberry, Prashant~K. Jha, J.~Tinsley Oden, and Omar Ghattas.
\newblock {Residual-based error correction for neural operator accelerated infinite-dimensional Bayesian inverse problems}.
\newblock \emph{Journal of Computational Physics}, 486:\penalty0 112104, 2023.
\newblock \doi{10.1016/j.jcp.2023.112104}.

\bibitem[Deuflhard(2004)]{deuflhard2004newton}
Peter Deuflhard.
\newblock \emph{Newton Methods for Nonlinear Problems: Affine Invariance and Adaptive Algorithms}, volume~35 of \emph{Springer Series in Computational Mathematics}.
\newblock Springer, Berlin, 2004.

\bibitem[Gonzalez et~al.(2024)Gonzalez, Warrington, Smith, and Linderman]{gonzalez2024towards}
Xavier Gonzalez, Andrew Warrington, Jimmy T.~H. Smith, and Scott~W. Linderman.
\newblock Towards scalable and stable parallelization of nonlinear {RNNs}.
\newblock In \emph{Advances in Neural Information Processing Systems}, volume~37, 2024.

\bibitem[Goswami et~al.(2023)Goswami, Bora, Yu, and Karniadakis]{goswami2023physics}
Somdatta Goswami, Aniruddha Bora, Yue Yu, and George~Em Karniadakis.
\newblock {Physics-informed deep neural operator networks}.
\newblock In \emph{Machine Learning in Modeling and Simulation}, pages 219--254. Springer, Cham, 2023.
\newblock \doi{10.1007/978-3-031-36644-4_6}.

\bibitem[Hairer et~al.(1993)Hairer, N{\o}rsett, and Wanner]{hairer1993solving}
Ernst Hairer, Syvert~P. N{\o}rsett, and Gerhard Wanner.
\newblock \emph{{Solving Ordinary Differential Equations I: Nonstiff Problems}}.
\newblock Springer, Berlin, 2nd edition, 1993.
\newblock \doi{10.1007/978-3-540-78862-1}.

\bibitem[Hillebrecht and Unger(2022)]{hillebrecht2022certified}
Birgit Hillebrecht and Benjamin Unger.
\newblock {Certified machine learning: a posteriori error estimation for physics-informed neural networks}.
\newblock In \emph{2022 International Joint Conference on Neural Networks (IJCNN)}, pages 1--8. IEEE, 2022.
\newblock \doi{10.1109/IJCNN55064.2022.9892569}.

\bibitem[Jha(2024)]{jha2024residual}
Prashant~K. Jha.
\newblock {Residual-based error corrector operator to enhance accuracy and reliability of neural operator surrogates of nonlinear variational boundary-value problems}.
\newblock \emph{Computer Methods in Applied Mechanics and Engineering}, 419:\penalty0 116595, 2024.
\newblock \doi{10.1016/j.cma.2023.116595}.

\bibitem[Kantorovich and Akilov(1982)]{kantorovich1982functional}
Leonid~V. Kantorovich and Gleb~P. Akilov.
\newblock \emph{Functional Analysis}.
\newblock Pergamon Press, Oxford, 2nd edition, 1982.

\bibitem[Kouri et~al.(2013)Kouri, Heinkenschloss, Ridzal, and van Bloemen~Waanders]{kouri2013trust}
D.~P. Kouri, M.~Heinkenschloss, D.~Ridzal, and B.~G. van Bloemen~Waanders.
\newblock {A Trust-Region Algorithm with Adaptive Stochastic Collocation for PDE Optimization under Uncertainty}.
\newblock \emph{SIAM Journal on Scientific Computing}, 35\penalty0 (4), 2013.
\newblock \doi{10.1137/120892362}.

\bibitem[Kurutach et~al.(2018)Kurutach, Clavera, Duan, Tamar, and Abbeel]{kurutach2018model}
Thanard Kurutach, Ignasi Clavera, Yan Duan, Aviv Tamar, and Pieter Abbeel.
\newblock {Model-ensemble trust-region policy optimization}.
\newblock In \emph{International Conference on Learning Representations (ICLR)}, 2018.
\newblock arXiv:1802.10592.

\bibitem[Lambrianides et~al.(2020)Lambrianides, Gong, and Venturi]{lambrianides2020scalable}
Panos Lambrianides, Qi~Gong, and Daniele Venturi.
\newblock A new scalable algorithm for computational optimal control under uncertainty.
\newblock \emph{Journal of Computational Physics}, 420:\penalty0 109710, 2020.
\newblock \doi{10.1016/j.jcp.2020.109710}.

\bibitem[Li et~al.(2024)Li, Zheng, Kovachki, Jin, Chen, Liu, Azizzadenesheli, and Anandkumar]{li2024physics}
Zongyi Li, Hongkai Zheng, Nikola Kovachki, David Jin, Haoxuan Chen, Burigede Liu, Kamyar Azizzadenesheli, and Anima Anandkumar.
\newblock {Physics-informed neural operator for learning partial differential equations}.
\newblock \emph{ACM/IMS Journal of Data Science}, 1\penalty0 (3):\penalty0 1--27, 2024.
\newblock \doi{10.1145/3648506}.

\bibitem[Lim et~al.(2024)Lim, Zhu, Selfridge, and Kasim]{lim2024parallelizing}
Yi~Heng Lim, Qi~Zhu, Joshua Selfridge, and Muhammad~Firmansyah Kasim.
\newblock {Parallelizing non-linear sequential models over the sequence length}.
\newblock In \emph{International Conference on Learning Representations (ICLR)}, 2024.
\newblock arXiv:2309.12252.

\bibitem[Lu et~al.(2021)Lu, Jin, Pang, Zhang, and Karniadakis]{lu2021deeponet}
Lu~Lu, Pengzhan Jin, Guofei Pang, Zhongqiang Zhang, and George~Em Karniadakis.
\newblock Learning nonlinear operators via {DeepONet} based on the universal approximation theorem of operators.
\newblock \emph{Nature Machine Intelligence}, 3\penalty0 (3):\penalty0 218--229, 2021.
\newblock \doi{10.1038/s42256-021-00302-5}.

\bibitem[Lundqvist and Oliveira(2025)]{lundqvist2025employing}
Oliver G.~S. Lundqvist and Fabricio Oliveira.
\newblock {Employing deep neural operators for PDE control by decoupling training and optimization}.
\newblock arXiv:2506.04742, 2025.

\bibitem[Luo et~al.(2025)Luo, O'Leary-Roseberry, Chen, and Ghattas]{luo2025efficient}
Dingcheng Luo, Thomas O'Leary-Roseberry, Peng Chen, and Omar Ghattas.
\newblock {Efficient PDE-Constrained Optimization Under High-Dimensional Uncertainty Using Derivative-Informed Neural Operators}.
\newblock \emph{SIAM Journal on Scientific Computing}, 47\penalty0 (4):\penalty0 C899--C931, 2025.
\newblock \doi{10.1137/23m157956x}.

\bibitem[Lye et~al.(2021)Lye, Mishra, Ray, and Chandrashekar]{lye2021iterative}
Kjetil~O. Lye, Siddhartha Mishra, Deep Ray, and Praveen Chandrashekar.
\newblock {Iterative surrogate model optimization (ISMO): An active learning algorithm for PDE constrained optimization with deep neural networks}.
\newblock \emph{Computer Methods in Applied Mechanics and Engineering}, 374:\penalty0 113575, 2021.
\newblock \doi{10.1016/j.cma.2020.113575}.

\bibitem[Mandl et~al.(2026)Mandl, Nayak, Ricken, and Goswami]{mandl2026physics}
Luis Mandl, Dibyajyoti Nayak, Tim Ricken, and Somdatta Goswami.
\newblock {Physics-informed time-integrated DeepONet: Temporal tangent space operator learning for high-accuracy inference}.
\newblock \emph{Computer Methods in Applied Mechanics and Engineering}, 455:\penalty0 118917, 2026.
\newblock \doi{10.1016/j.cma.2026.118917}.

\bibitem[McCabe et~al.(2023)McCabe, Harrington, Subramanian, and Brown]{mccabe2023towards}
Michael McCabe, Peter Harrington, Shashank Subramanian, and Jed Brown.
\newblock {Towards stability of autoregressive neural operators}.
\newblock \emph{Transactions on Machine Learning Research}, 2023.
\newblock arXiv:2306.10619.

\bibitem[Melnikov and Milz(2026)]{melnikov2026convergence}
Olena Melnikov and Johannes Milz.
\newblock {Convergence Rates for Ensemble-Based Solutions to Optimal Control of Uncertain Dynamical Systems}.
\newblock \emph{SIAM Journal on Control and Optimization}, 64\penalty0 (3):\penalty0 1631--1650, 2026.
\newblock \doi{10.1137/24m1693076}.

\bibitem[O'Leary-Roseberry et~al.(2024)O'Leary-Roseberry, Chen, Villa, and Ghattas]{oleary2024derivative}
Thomas O'Leary-Roseberry, Peng Chen, Umberto Villa, and Omar Ghattas.
\newblock {Derivative-Informed Neural Operator: An efficient framework for high-dimensional parametric derivative learning}.
\newblock \emph{Journal of Computational Physics}, 496:\penalty0 112555, 2024.
\newblock \doi{10.1016/j.jcp.2023.112555}.

\bibitem[Phelps et~al.(2014)Phelps, Gong, Royset, Walton, and Kaminer]{phelps2014consistent}
Chris Phelps, Qi~Gong, Johannes~O. Royset, Claire Walton, and Isaac Kaminer.
\newblock {Consistent approximation of a nonlinear optimal control problem with uncertain parameters}.
\newblock \emph{Automatica}, 50\penalty0 (12):\penalty0 2987--2997, 2014.
\newblock \doi{10.1016/j.automatica.2014.10.025}.

\bibitem[Phelps et~al.(2016)Phelps, Royset, and Gong]{phelps2016optimal}
Chris Phelps, Johannes~O. Royset, and Qi~Gong.
\newblock {Optimal Control of Uncertain Systems Using Sample Average Approximations}.
\newblock \emph{SIAM Journal on Control and Optimization}, 54\penalty0 (1):\penalty0 1--29, 2016.
\newblock \doi{10.1137/140983161}.

\bibitem[Qian et~al.(2017)Qian, Grepl, Veroy, and Willcox]{qian2017certified}
Elizabeth Qian, Martin Grepl, Karen Veroy, and Karen Willcox.
\newblock {A Certified Trust Region Reduced Basis Approach to PDE-Constrained Optimization}.
\newblock \emph{SIAM Journal on Scientific Computing}, 39\penalty0 (5):\penalty0 S434--S460, 2017.
\newblock \doi{10.1137/16m1081981}.

\bibitem[Rockafellar and Uryasev(2000)]{rockafellar2000optimization}
R.~Tyrrell Rockafellar and Stanislav Uryasev.
\newblock {Optimization of conditional value-at-risk}.
\newblock \emph{The Journal of Risk}, 2\penalty0 (3):\penalty0 21--41, 2000.
\newblock \doi{10.21314/jor.2000.038}.

\bibitem[Ross et~al.(2015)Ross, Proulx, Karpenko, and Gong]{ross2015riemann}
I.~Michael Ross, Ronald~J. Proulx, Mark Karpenko, and Qi~Gong.
\newblock {Riemann--Stieltjes Optimal Control Problems for Uncertain Dynamic Systems}.
\newblock \emph{Journal of Guidance, Control, and Dynamics}, 38\penalty0 (7):\penalty0 1251--1263, 2015.
\newblock \doi{10.2514/1.g000505}.

\bibitem[S{\"a}rkk{\"a} and Garc{\'i}a-Fern{\'a}ndez(2023)]{sarkka2023temporal}
Simo S{\"a}rkk{\"a} and {\'A}ngel~F. Garc{\'i}a-Fern{\'a}ndez.
\newblock Temporal parallelization of dynamic programming and linear quadratic control.
\newblock \emph{IEEE Transactions on Automatic Control}, 68\penalty0 (2):\penalty0 851--866, 2023.

\bibitem[Shaffer et~al.(2016)Shaffer, Karpenko, and Gong]{shaffer2016unscented}
Richard Shaffer, Mark Karpenko, and Qi~Gong.
\newblock {Unscented guidance for waypoint navigation of a fixed-wing UAV}.
\newblock In \emph{2016 American Control Conference (ACC)}. IEEE, 2016.
\newblock \doi{10.1109/acc.2016.7524959}.

\bibitem[Shapiro et~al.(2009)Shapiro, Dentcheva, and Ruszczy{\'n}ski]{shapiro2009lectures}
Alexander Shapiro, Darinka Dentcheva, and Andrzej Ruszczy{\'n}ski.
\newblock \emph{Lectures on Stochastic Programming: Modeling and Theory}.
\newblock MPS-SIAM Series on Optimization. Society for Industrial and Applied Mathematics, Philadelphia, 2009.
\newblock \doi{10.1137/1.9780898718751}.

\bibitem[Skeel(1986)]{skeel1986thirteen}
Robert~D. Skeel.
\newblock {Thirteen ways to estimate global error}.
\newblock \emph{Numerische Mathematik}, 48\penalty0 (1):\penalty0 1--20, 1986.
\newblock \doi{10.1007/BF01389440}.

\bibitem[Stetter(1978)]{stetter1978defect}
Hans~J. Stetter.
\newblock {The defect correction principle and discretization methods}.
\newblock \emph{Numerische Mathematik}, 29\penalty0 (4):\penalty0 425--443, 1978.
\newblock \doi{10.1007/BF01432879}.

\bibitem[Tancik et~al.(2020)Tancik, Srinivasan, Mildenhall, Fridovich-Keil, Raghavan, Singhal, Ramamoorthi, Barron, and Ng]{tancik2020fourier}
Matthew Tancik, Pratul~P. Srinivasan, Ben Mildenhall, Sara Fridovich-Keil, Nithin Raghavan, Utkarsh Singhal, Ravi Ramamoorthi, Jonathan~T. Barron, and Ren Ng.
\newblock {Fourier features let networks learn high frequency functions in low dimensional domains}.
\newblock In \emph{Advances in Neural Information Processing Systems}, volume~33, pages 7537--7547, 2020.

\bibitem[W{\"a}chter and Biegler(2006)]{wachter2006implementation}
Andreas W{\"a}chter and Lorenz~T. Biegler.
\newblock {On the implementation of an interior-point filter line-search algorithm for large-scale nonlinear programming}.
\newblock \emph{Mathematical Programming}, 106\penalty0 (1):\penalty0 25--57, 2006.
\newblock \doi{10.1007/s10107-004-0559-y}.

\bibitem[Wang and Perdikaris(2023)]{wang2023long}
Sifan Wang and Paris Perdikaris.
\newblock {Long-time integration of parametric evolution equations with physics-informed DeepONets}.
\newblock \emph{Journal of Computational Physics}, 475:\penalty0 111855, 2023.
\newblock \doi{10.1016/j.jcp.2022.111855}.

\bibitem[Wang et~al.(2021{\natexlab{a}})Wang, Bhouri, and Perdikaris]{wang2021fast}
Sifan Wang, Mohamed~Aziz Bhouri, and Paris Perdikaris.
\newblock Fast {PDE}-constrained optimization via self-supervised operator learning.
\newblock \emph{arXiv preprint arXiv:2110.13297}, 2021{\natexlab{a}}.

\bibitem[Wang et~al.(2021{\natexlab{b}})Wang, Teng, and Perdikaris]{wang2021understanding}
Sifan Wang, Yujun Teng, and Paris Perdikaris.
\newblock {Understanding and mitigating gradient flow pathologies in physics-informed neural networks}.
\newblock \emph{SIAM Journal on Scientific Computing}, 43\penalty0 (5):\penalty0 A3055--A3081, 2021{\natexlab{b}}.
\newblock \doi{10.1137/20M1318043}.

\bibitem[Wang et~al.(2021{\natexlab{c}})Wang, Wang, and Perdikaris]{wang2021learning}
Sifan Wang, Hanwen Wang, and Paris Perdikaris.
\newblock {Learning the solution operator of parametric partial differential equations with physics-informed DeepONets}.
\newblock \emph{Science Advances}, 7\penalty0 (40):\penalty0 eabi8605, 2021{\natexlab{c}}.
\newblock \doi{10.1126/sciadv.abi8605}.

\bibitem[Wang et~al.(2024)Wang, Sankaran, and Perdikaris]{wang2024respecting}
Sifan Wang, Shyam Sankaran, and Paris Perdikaris.
\newblock {Respecting causality for training physics-informed neural networks}.
\newblock \emph{Computer Methods in Applied Mechanics and Engineering}, 421:\penalty0 116813, 2024.
\newblock \doi{10.1016/j.cma.2024.116813}.

\bibitem[Yao et~al.(2025)Yao, Luo, Cao, Kovachki, O'Leary-Roseberry, and Ghattas]{yao2025derivative}
Boyuan Yao, Dingcheng Luo, Lianghao Cao, Nikola Kovachki, Thomas O'Leary-Roseberry, and Omar Ghattas.
\newblock {Derivative-informed Fourier neural operator: universal approximation and applications to PDE-constrained optimization}.
\newblock arXiv:2512.14086, 2025.

\bibitem[Zadunaisky(1976)]{zadunaisky1976estimation}
Pedro~E. Zadunaisky.
\newblock {On the estimation of errors propagated in the numerical integration of ordinary differential equations}.
\newblock \emph{Numerische Mathematik}, 27\penalty0 (1):\penalty0 21--39, 1976.
\newblock \doi{10.1007/BF01399082}.

\bibitem[Zhu et~al.(2023)Zhu, Zhang, Jiao, Karniadakis, and Lu]{zhu2023reliable}
Min Zhu, Handi Zhang, Anran Jiao, George~Em Karniadakis, and Lu~Lu.
\newblock {Reliable extrapolation of deep neural operators informed by physics or sparse observations}.
\newblock \emph{Computer Methods in Applied Mechanics and Engineering}, 412:\penalty0 116064, 2023.
\newblock \doi{10.1016/j.cma.2023.116064}.

\end{thebibliography}
\clearpage
\appendix

\section{Definitions, proofs and derivations}\label{app:proofs}

\subsection*{Propagated curvature and the certificate constants}
For $0 \le \tau \le t \le t_f$ and the convex sets $\mathcal C(\tau)$ of Lemma~\ref{lem:identity} let
\begin{equation}\label{eq:K}
\begin{aligned}
  \hat k(t, \tau) &:= \sup_{\zeta \in \mathcal C(\tau)} \, \sup_{\|h\| \le 1}
  \Big\| \hat\Phi(t, \tau) \Big( h^\top \partial_s^2 f_i\big(\zeta, u(\tau), \xi\big)\, h \Big)_{i=1}^{n_s} \Big\|, \qquad
  \hat K(t) := \int_0^t \hat k(t, \tau)\, d\tau .
\end{aligned}
\end{equation}
Here $h^\top \partial_s^2 f_i\, h$ is the curvature of the $i$-th component of $f$ in the direction $h$. The quantity $\hat k(t, \tau)$ takes this curvature in the worst direction $h$ and at the worst state $\zeta$ in $\mathcal C(\tau)$, a set that contains the prediction, the true state, and all states between them, and carries it from time $\tau$ to time $t$ with the transition matrix $\hat\Phi(t, \tau)$ of the linearized dynamics. The constant $\hat K(t)$ adds these propagated curvatures over all earlier times $\tau \le t$.

Similarly, along the corrected state, let $\bar\Phi$ solve $\partial_t \bar\Phi(t, \tau) = \bar A(t)\, \bar\Phi(t, \tau)$ with $\bar\Phi(\tau, \tau) = I$, let $\bar{\mathcal C}(\tau) := \{\zeta \in \mathbb{R}^{n_s} \colon \|\zeta - \bar s(\tau)\| \le 2\bar a\}$ with $\bar a$ of \eqref{eq:abar}, and
\begin{equation}\label{eq:K_bar}
  \bar k(t, \tau) := \sup_{\zeta \in \bar{\mathcal C}(\tau)} \, \sup_{\|h\| \le 1}
  \Big\| \bar\Phi(t, \tau) \Big( h^\top \partial_s^2 f_i\big(\zeta, u(\tau), \xi\big)\, h \Big)_{i=1}^{n_s} \Big\|, \qquad
  \bar K := \sup_{t \le t_f} \int_0^t \bar k(t, \tau)\, d\tau .
\end{equation}

\subsection*{Proofs}

\renewcommand{\proofname}{Proof of Lemma~\ref{lem:identity}}\begin{proof}

Subtracting \eqref{eq:estimate_ode} from \eqref{eq:error_ode} gives
\begin{equation}\label{eq:true_error_diff}
    \tfrac{d}{dt}(e - \hat e) = \hat A\,(e - \hat e) - q,
\end{equation}
with $(e - \hat e)(0) = 0$. Writing $\bar s - s = e - \hat e$ and using the state transition matrix to solve \eqref{eq:true_error_diff}, we get
\begin{equation}\label{eq:identity}
  \bar s(t) - s(t) = -\int_0^t \hat\Phi(t, \tau)\, q(\tau)\, d\tau .
\end{equation}
For every $\tau \le t$ and $\theta \in [0, 1]$, the point $\hat s(\tau) - \theta e(\tau)$ lies on the segment between $\hat s(\tau)$ and $s(\tau)$, and hence in the convex set $\mathcal C(\tau)$. With $h = e(\tau)/\|e(\tau)\|$ in \eqref{eq:K} and $\int_0^1 (1 - \theta)\, d\theta = \tfrac12$ we obtain
\begin{equation}\label{eq:remainder_bound}
\begin{aligned}
  \big\| \hat\Phi(t, \tau)\, q(\tau) \big\|
  &\le \int_0^1 (1 - \theta)\, \Big\| \hat\Phi(t, \tau) \Big( e(\tau)^\top \partial_s^2 f_i\big(\hat s(\tau) - \theta e(\tau), u(\tau), \xi\big)\, e(\tau) \Big)_{i=1}^{n_s} \Big\|\, d\theta \\
  &\le \tfrac12\, \hat k(t, \tau)\, \|e(\tau)\|^2 .
\end{aligned}
\end{equation}
Taking the norm of \eqref{eq:identity}, the triangle inequality for integrals and \eqref{eq:remainder_bound} give
\begin{equation}\label{eq:correction_steps}
\begin{aligned}
  \big\| \bar s(t) - s(t) \big\|
  &\le \int_0^t \big\| \hat\Phi(t, \tau)\, q(\tau) \big\|\, d\tau
  \le \tfrac12 \int_0^t \hat k(t, \tau)\, \|e(\tau)\|^2\, d\tau \\
  &\le \tfrac12\, \sup_{\tau \le t} \|e(\tau)\|^2 \int_0^t \hat k(t, \tau)\, d\tau
  = \tfrac12\, \hat K(t) \sup_{\tau \le t} \|e(\tau)\|^2 ,
\end{aligned}
\end{equation}
which is \eqref{eq:correction}. For \eqref{eq:error_bound}, the solution formula for \eqref{eq:error_ode} with $e(0) = 0$ is
\[
  e(t) = \int_0^t \hat\Phi(t, \tau)\, \big[ R(\tau) - q(\tau) \big]\, d\tau
  = \hat e(t) - \int_0^t \hat\Phi(t, \tau)\, q(\tau)\, d\tau ,
\]
where the residual part is $\hat e(t)$ by \eqref{eq:estimate}. The norm of the remainder part is bounded by \eqref{eq:correction_steps}, so
\[
  \|e(t)\| \le \|\hat e(t)\| + \tfrac12\, \hat K(t) \sup_{\tau \le t} \|e(\tau)\|^2 .
\]
\end{proof}

\renewcommand{\proofname}{Proof of Theorem~\ref{thm:certificate}}\begin{proof}

Since $\hat e(0) = 0$, the corrected state starts at $\bar s(0) = s_0(\xi)$, and $\dot{\bar s} = f(\hat s) + R - \hat A \hat e - R$ by \eqref{eq:residual} and \eqref{eq:estimate_ode}, so $\bar R = \dot{\bar s} - f(\bar s) = f(\hat s) - \hat A\, \hat e - f(\bar s)$ is its residual. The proof of Lemma~\ref{lem:identity} uses only the start state, the residual and Taylor's theorem on the sets, so the lemma holds for $\bar s$ in place of $\hat s$, with $\bar R$, $\bar\Phi$, $\bar k$ and sets that contain $s(\tau)$ and $\bar s(\tau)$. Let $\bar\eta(t) := \sup_{\tau \le t} \|\bar s(\tau) - s(\tau)\|$. On every interval $[0, t]$ on which $s(\tau) \in \bar{\mathcal C}(\tau)$, the inequality \eqref{eq:error_bound} for $\bar s$ gives
\begin{equation}\label{eq:quadratic}
  \bar\eta(t) \le \bar a + \tfrac12\, \bar K\, \bar\eta(t)^2 ,
\end{equation}
since its first-order term is at most $\bar a$ by \eqref{eq:abar} and its right-hand side is nondecreasing in $t$. Under \eqref{eq:certificate} the polynomial $\tfrac12 \bar K \eta^2 - \eta + \bar a$ has the real roots $\bar\eta_\pm = \big(1 \pm \sqrt{1 - 2 \bar K \bar a}\big)/\bar K$, and $\bar\eta_- < 2 \bar a$ since $1 - \sqrt{1 - x} < x$ for $0 < x < 1$. The polynomial $\tfrac12 \bar K \eta^2 - \eta + \bar a$ is negative strictly between its roots, so \eqref{eq:quadratic} excludes the values $\bar\eta_- < \bar\eta < \bar\eta_+$. Since $\bar\eta$ is continuous with $\bar\eta(0) = 0 < \bar\eta_-$, it cannot cross from below $\bar\eta_-$ to above $\bar\eta_+$, so $\bar\eta \le \bar\eta_- < 2\bar a$ as long as \eqref{eq:quadratic} holds. The strict inequality keeps $s$ away from the boundary of the tubes, so $s$ continues to $t_f$ and \eqref{eq:quadratic} holds on all of $[0, t_f]$, which closes the continuation.
\end{proof}

\renewcommand{\proofname}{Proof of Theorem~\ref{thm:honest}}\begin{proof}

Adding and subtracting $J_M(u)$ and using the triangle inequality,
\begin{equation*}
  \big| J(u) - \bar J_M(u) \big| \le \big| J(u) - J_M(u) \big| + \big| J_M(u) - \bar J_M(u) \big| .
\end{equation*}
The costs $J_M$ and $\bar J_M$ share the control cost $c(u)$ and differ only in the end states, so with the Lipschitz continuity of $\ell$
\begin{equation*}
  \big| J_M(u) - \bar J_M(u) \big|
  = \Big| \frac{1}{M} \sum_{i=1}^M \Big( \ell\big(s_i(t_f; u)\big) - \ell\big(\bar s_i(t_f; u)\big) \Big) \Big|
  \le \frac{L_\ell}{M} \sum_{i=1}^M \big\| \bar s_i(t_f; u) - s_i(t_f; u) \big\| ,
\end{equation*}
which gives \eqref{eq:honest}, and Theorem~\ref{thm:certificate} bounds the summands by $2\bar a_i$ for the certified scenarios. For \eqref{eq:subopt}, adding and subtracting $\bar J_M(\hat u)$ and $\bar J_M(u^*)$ gives
\begin{equation*}
  J(\hat u) - J(u^*)
  = \big[ J(\hat u) - \bar J_M(\hat u) \big] + \big[ \bar J_M(\hat u) - \bar J_M(u^*) \big] + \big[ \bar J_M(u^*) - J(u^*) \big] .
\end{equation*}
The middle bracket is at most $\varepsilon_{\mathrm{opt}}$, since $\hat u$ minimizes $\bar J_M$ over $\mathcal U_{\mathcal W}$ up to $\varepsilon_{\mathrm{opt}}$ and $u^* \in \mathcal U_{\mathcal W}$. The outer brackets are bounded by \eqref{eq:honest} at $\hat u$ and at $u^*$, with $|J - J_M|$ at both controls bounded by its supremum over $\mathcal U_{\mathcal W}$, which gives the factor 2 and the two correction terms of \eqref{eq:subopt}.
\end{proof}

\subsection*{Newton form of the correction}
Let $G(s) := \dot s - f(s)$ on trajectories with $s(0) = s_0$, so that the true state is the root of $G$ and $G(\hat s) = R$. The derivative of $G$ at $\hat s$ is the linear operator
\begin{equation*}
  G'(\hat s)[\delta] = \dot\delta - \hat A\,\delta, \qquad \delta(0) = 0,
\end{equation*}
since $s \mapsto \dot s$ is linear and the derivative of $s \mapsto f(s)$ at $\hat s$ is $\hat A$. One Newton step from $\hat s$ solves $G(\hat s) + G'(\hat s)[\,s^1 - \hat s\,] = 0$, that is
\begin{equation*}
  R + \tfrac{d}{dt}(s^1 - \hat s) - \hat A\,(s^1 - \hat s) = 0, \qquad (s^1 - \hat s)(0) = 0 .
\end{equation*}
With $s^1 - \hat s = -\hat e$, this is the estimate equation \eqref{eq:estimate_ode}, and a linear initial value problem has one solution, so $s^1 = \hat s - \hat e = \bar s$. The correction is one Newton step on the dynamics, and Lemma~\ref{lem:identity} is the quadratic convergence of that step, with the propagated curvature $\hat K$ as its constant. The residual of $\bar s$ is $\bar R = \dot{\bar s} - f(\bar s) = f(\hat s) - \hat A\,\hat e - f(\bar s)$, which Taylor's theorem at $\hat s$ makes of second order in $\hat e$, so the next step, $-\int_0^t \bar\Phi(t, \tau)\, \bar R(\tau)\, d\tau$, is bounded by $\bar a$ of \eqref{eq:abar}, and condition \eqref{eq:certificate} is a Newton--Kantorovich type condition at $\bar s$ \citep{deuflhard2004newton}.

\clearpage

\section{Problem formulation of the solvers}\label{app:problem}

\subsection*{Sample-average problem}
All solvers approximate the same problem. With $M$ scenarios $\xi_1, \dots, \xi_M$, the sample-average approximation of \eqref{eq:ocp} for the UAV problem reads
\begin{equation}\label{eq:uav_saa}
\begin{aligned}
  \min_{u \in \mathcal U} \quad & J_M(u) = \frac{1}{M} \sum_{i=1}^{M} \big\| p_i(t_f) - p^* \big\|^2
    + \lambda_u \int_0^{t_f} \|\dot u(t)\|^2 \, dt \\
  \text{s.t.} \quad & \dot s_i(t) = f\big(s_i(t), u(t), \xi_i\big), \qquad s_i(0) = s_0(\xi_i), \qquad i = 1, \dots, M, \\
  & u(0) = u_0, \qquad 3\ \mathrm{N} \le T(t) \le 35\ \mathrm{N}, \qquad |\alpha(t)| \le \pi/12, \\
  & |\dot T(t)| \le 1\ \mathrm{N/s}, \qquad |\dot\alpha(t)| \le 0.05\ \mathrm{rad/s}, \qquad |\dot\mu(t)| \le 0.05\ \mathrm{rad/s}, \\
  & v_{\min} \le \frac{1}{M} \sum_{i=1}^{M} v_i(t) \le v_{\max}, \qquad
    \Big| \frac{1}{M} \sum_{i=1}^{M} \gamma_i(t) \Big| \le \gamma_{\max}, \qquad t \in [0, t_f].
\end{aligned}
\end{equation}
Here $s_i = (p_i, v_i, \gamma_i, \sigma_i)$ is the state of scenario $\xi_i \in \mathbb{R}^{10}$ under the flight equations \eqref{eq:uav}, with the start state $s_0(\xi_i) = (p_{0,i}, v_{0,i}, \gamma_{0,i}, \sigma_{0,i})$ and the aerodynamic coefficients $C_{x0,i}$, $C_{x\alpha,i}$, $C_{z0,i}$ and $C_{z\alpha,i}$ of $\xi_i$ in \eqref{eq:aero}. The scenarios are drawn independently with the ranges of Table~\ref{tab:setup}. The input $u = (T, \alpha, \mu)$ is the same for every scenario, and Table~\ref{tab:setup} gives $u_0$, $v_{\min}$, $v_{\max}$ and $\gamma_{\max}$.

\subsection*{Problem parameters and uncertainty domains}
The parameters and uncertainty set for problem \eqref{eq:uav_saa} are given in Table~\ref{tab:setup}.
\begin{table}[htbp]\centering\footnotesize
\caption{Problem data of the UAV benchmark.}\label{tab:setup}
\renewcommand{\arraystretch}{1.05}
\begin{tabular}{@{}l>{\raggedright\arraybackslash}p{0.84\textwidth}@{}}
\toprule
constants & $m_a = 2$ kg,\; $S = 0.982$ m$^2$,\; $g_0 = 9.8$ m/s$^2$,\; $\varrho(z) = 1.21\, e^{-z/8000}$ kg/m$^3$ \\
\midrule
task & $p^* = (500, 500, 500)$ m,\; $t_f = 50$ s,\; $\lambda_u = 1/200$ \\
\midrule
scenario $\xi$ & $\|p(0)\| \sim U(0, 5)$ m with uniform spherical angles,\; $v_0 \sim U(25.575, 29.425)$ m/s,\; $\gamma_0 \sim U(-0.05, 0.05)$,\; $\sigma_0 \sim U(3.1, 3.2)$,\; $C_{x0} \sim U(-0.0380, -0.0330)$,\; $C_{x\alpha} \sim U(0.0027, 0.0031)$,\; $C_{z0} \sim U(-0.0589, -0.0548)$,\; $C_{z\alpha} \sim U(-5.9685, -5.1875)$ \\
\midrule
admissible set $\mathcal U$ & $u(0) = u_0 = (16.1\ \mathrm{N}, -0.0088, 0)$,\; $3 \le T \le 35$ N,\; $|\alpha| \le \pi/12$,\; $|\dot T| \le 1$ N/s,\; $|\dot\alpha|, |\dot\mu| \le 0.05$ rad/s \\
\midrule
mean state & $v_{\min} = 13$ m/s,\; $v_{\max} = 42$ m/s,\; $\gamma_{\max} = \pi/6$ \\
\bottomrule
\end{tabular}
\end{table}

\subsection*{Spline inputs}
All solvers use the spline inputs of Section~\ref{sec:single_shooting},
\begin{equation*}\tag{\ref{eq:spline}}
  u(t; w) = u_{\mathrm{mid}} + u_{\mathrm{half}} \odot \sum_{k=1}^{n_b} w_k\, B_k(t),
\end{equation*}
with $n_b = 26$ cubic B-splines per input on uniform knots 2.2\,s apart, the coefficients $w \in \mathcal W = [-1, 1]^{3 \times 26}$, and the centre $u_{\mathrm{mid}}$ and the half-width $u_{\mathrm{half}}$ of Table~\ref{tab:spline_ranges}. The first coefficient of each input is fixed so that $u(0) = u_0$, which leaves 75 free coefficients. Since the B-splines are nonnegative and sum to one, $u(t; w)$ stays within the ranges of Table~\ref{tab:spline_ranges} for every $w \in \mathcal W$, inside the bounds on $T$ and $\alpha$. The benchmark has no bound on the bank angle, and the spline limits it to $|\mu| \le 0.8$\,rad, the range of the training controls. The inputs are evaluated on a grid of step 0.1\,s, and the input rates are the difference quotients of the inputs on this grid.

\begin{table}[htbp]\centering\footnotesize
\caption{Centre, half-width and range of the spline inputs \eqref{eq:spline}.}\label{tab:spline_ranges}
\begin{tabular}{lccc}
\toprule
input & $u_{\mathrm{mid}}$ & $u_{\mathrm{half}}$ & range \\
\midrule
thrust $T$ & 19 N & $0.98 \times 16$ N & $[3.3, 34.7]$ N \\
angle of attack $\alpha$ & 0 & $0.98 \times \pi/12$ & $[-0.257, 0.257]$ rad \\
bank angle $\mu$ & 0 & 0.8 rad & $[-0.8, 0.8]$ rad \\
\bottomrule
\end{tabular}
\end{table}

\subsection*{State models and margin}
All solvers eliminate the dynamics via single-shooting. For $w \in \mathcal W$, the spline gives the input $u(\cdot\,; w)$, and a state model gives the states of all $M$ scenarios. The operator solver uses the corrected states \eqref{eq:corrected} with the corrected cost $\bar J_M$ of \eqref{eq:checkedcost}, and the reference solvers use RK4 states with the cost $J_M$ of \eqref{eq:uav_saa}. Both state models evaluate the states on grids coarser than the held-out scenario evaluation at step 0.1\,s, and these grids underestimate the early peak of the mean flight-path angle by up to 0.03\,rad. Every solver therefore tightens the last bound of \eqref{eq:uav_saa} by a margin $\Delta_\gamma$ to
\begin{equation}\label{eq:margin}
  \Big| \frac{1}{M} \sum_{i=1}^{M} \gamma_i(t) \Big| \le \gamma_{\max} - \Delta_\gamma, \qquad \Delta_\gamma = 0.035\ \mathrm{rad}, \qquad t \in [0, t_f].
\end{equation}
This is the smallest margin at which every solve of Section~\ref{sec:exp} satisfies the bound of \eqref{eq:uav_saa} under RK4 at step 0.1\,s (Appendix~\ref{app:rk4steps}).

\subsection*{Reduced problem}
With the states eliminated, every bound of \eqref{eq:uav_saa}, with \eqref{eq:margin} in place of the last one, applies to a scalar function of time $y(\cdot\,; w)$ that depends on $w$ alone. We write such a bound as $y_{\min} \le y \le y_{\max}$ and collect the bounded inputs and mean states in $\mathcal Y = \{T,\, \alpha,\, \mathbb{E}_M[v],\, \mathbb{E}_M[\gamma]\}$, with the sample mean $\mathbb{E}_M[v] := \frac{1}{M} \sum_{i=1}^{M} v_i$, and the bounded input rates in $\mathcal Y' = \{\dot T,\, \dot\alpha,\, \dot\mu\}$, with the mean states from the state model of the solver. The reduced problem is
\begin{equation}\label{eq:reduced}
  \min_{w \in \mathcal W} \; J_M\big(u(\cdot\,; w)\big) \qquad \text{s.t.} \qquad y_{\min} \le y(t_k; w) \le y_{\max}, \qquad y \in \mathcal Y \cup \mathcal Y',
\end{equation}
at the nodes $t_k$ of the grid of each function. Ipopt solves \eqref{eq:reduced} with RK4 states and hard constraints.

\subsection*{Penalized problem}
The operator solver and the reference Adam solver replace the bounds of \eqref{eq:reduced} with barrier and penalty terms. Only the box $\mathcal W$ remains, and Adam keeps it by projection after every step. The operator solver solves
\begin{equation}\label{eq:penalized}
\begin{aligned}
  \min_{w \in \mathcal W} \quad \bar J_M\big(u(\cdot\,; w)\big)
  &+ \sum_{y \in \mathcal Y} \Big( 10 \operatorname{bar}\big(y(\cdot\,; w)\big) + 10^8 \operatorname{pen}\big(y(\cdot\,; w)\big) \Big)
  + 10^8 \sum_{y \in \mathcal Y'} \operatorname{pen}\big(y(\cdot\,; w)\big) \\
  &+ 10 \sum_{i=1}^{M} \Big( \operatorname{bar}\big(\bar v_i(\cdot\,; w)\big) + \operatorname{bar}\big(\bar\gamma_i(\cdot\,; w)\big) \Big),
\end{aligned}
\end{equation}
with the mean states in $\mathcal Y$ taken from the corrected states. The last sum keeps every corrected scenario in the range in which the operator was trained, with its speed 0.5\,m/s inside the speed bounds and $|\bar\gamma_i| \le 0.9$. For a function $y$ at the $N_y$ nodes $t_k$ of its grid, let $\eta_k = \big(y(t_k; w) - y_{\min}\big)/(y_{\max} - y_{\min})$. The barrier and the penalty are
\begin{equation*}
\begin{aligned}
  \operatorname{bar}(y) &= \frac{1}{N_y} \sum_{k=1}^{N_y} \big( r(\eta_k) + r(1 - \eta_k) \big), \\
  \operatorname{pen}(y) &= \frac{1}{N_y} \sum_{k=1}^{N_y} \big( \max(-\eta_k, 0)^2 + \max(\eta_k - 1, 0)^2 \big),
\end{aligned}
\end{equation*}
with the relaxed log barrier
\begin{equation*}
  r(\eta) =
  \begin{cases}
    -\log \eta, & \eta \ge \delta, \\[2pt]
    -\log \delta + \big(1 - \eta/\delta\big) + \tfrac12 \big((\eta - \delta)/\delta\big)^2, & \eta < \delta,
  \end{cases}
  \qquad \delta = 0.02 .
\end{equation*}
The barrier is $C^1$ and keeps a restoring gradient after a bound is crossed. The penalty enforces the bounds where the barrier alone lets the mean flight-path angle drift past its bound. The rate penalties act at 98\,\% of the rate bounds. The reference Adam solver solves \eqref{eq:penalized} with $J_M$ and the RK4 states in place of $\bar J_M$ and the corrected states. It replaces the last sum by a quadratic penalty on $|\gamma_i| > 1$, which keeps the integration away from the singularity of \eqref{eq:uav} at $\cos\gamma = 0$.
\clearpage

\section{Operator and reference solvers}\label{app:details}

\subsection*{Architecture}
The operator is the DeepONet \eqref{eq:operator} with $n_\psi = 256$ terms per state. The two subnetworks, branch and trunk, are the modified multilayer perceptrons (MLPs) with skip connections described by \citet{wang2021understanding}. The branch network has four hidden layers of width 512 and sees the three input profiles sampled every 0.5 s, i.e., 303 values. The trunk has three hidden layers of width 512. The frame is turned to the start position and heading of every scenario, so $x_0$, $y_0$ and $\sigma_0$ do not enter the trunk, and the scenario variables reduce to seven. The trunk input is the scaled time $t/t_f$ with 48 Fourier feature pairs of it \citep{tancik2020fourier}, i.e., 96 values, the three start states $z_0$, $v_0$ and $\gamma_0$, and the four aerodynamic coefficients $C_{x0}$, $C_{x\alpha}$, $C_{z0}$ and $C_{z\alpha}$, 104 values in all. The output of the dot product of branch and trunk is the six states $(x,y,z,v,\gamma,\sigma)$. The state scales of \eqref{eq:operator} are $d = (1000, 1000, 500, 10, 0.5, 3)$. The architectural parameters are summarized in Table~\ref{tab:architecture}.

\begin{table}[htbp]\centering\footnotesize
\caption{Architecture of the operator. W $\times$ D is the hidden width and the number of hidden layers, and each state has a block of $n_\psi$ outputs.}
\label{tab:architecture}
\begin{tabular}{cccccccc}
\toprule
branch input & branch (W $\times$ D) & trunk input & trunk (W $\times$ D) & $n_\psi$ & outputs & activation & parameters \\
\midrule
303 & $512 \times 4$ & 104 & $512 \times 3$ & 256 & $6 \times 256$ & tanh & 3.52 M \\
\bottomrule
\end{tabular}
\end{table}

\subsection*{Training input}
The operator is trained on input profiles alone, with no trajectories and no task. The profiles come from a randomized autopilot, a proportional controller on the input rates that tracks random reference signals with random gains and random input noise. Every profile is one flight of $t_f = 50$ s on a scenario from the training box, computed at the step $\Delta t = 0.1$ s, and only the inputs of the flight are kept, i.e., the control. We do this to obtain controls that are feasible, i.e., flyable, and respect the bounds of \citet{lambrianides2020scalable}.

We obtain the feasible controls used in training by the following simple computational scheme. We create the feasible controls once in advance before the training iterations. Let $g_\ell$ denote a zero-mean Gaussian random field on $[0, t_f]$ with length scale $\ell$. Each flight draws three reference signals
\begin{equation*}
  v_r(t) = c_v + a_v\, g_{\ell_v}(t), \qquad
  \gamma_r(t) = c_\gamma + a_\gamma\, g_{\ell_\gamma}(t), \qquad
  \mu_r(t) = c_\mu + a_\mu\, g_{\ell_\mu}(t),
\end{equation*}
each clipped to its range, and three noise signals $n_j(t) = b_j\, g_{\ell_j}(t)$ for $j \in \{T, \alpha, \mu\}$. The constants, amplitudes, length scales and gains are uniform in the ranges of Table~\ref{tab:controls}, drawn once per flight.

\begin{table}[htbp]\centering\footnotesize
\caption{Generation of the training controls. Every range is sampled uniformly, once per flight.}
\label{tab:controls}
\setlength{\tabcolsep}{2pt}
\begin{tabular}{llll}
\toprule
reference & constant $c$ & amplitude $a$ & clipped to \\
\midrule
speed $v_r$ [m/s] & $[18, 36]$ & $[0, 10]$ & $[15, 40]$ \\
flight-path angle $\gamma_r$ [rad] & $[-0.3, 0.3]$ & $[0, 0.3]$ & $[-0.5, 0.5]$ \\
bank angle $\mu_r$ [rad] & $[-0.5, 0.5]$ & $[0, 0.5]$ & $[-0.8, 0.8]$ \\
\midrule
input & noise amplitude $b$ & rate bound $r$ & input box \\
\midrule
thrust $T$ & $[0, 0.4]$ N/s & 1 N/s & $[3, 35]$ N \\
angle of attack $\alpha$ & $[0, 0.01]$ rad/s & 0.05 rad/s & $[-\pi/8, \pi/8]$ rad \\
bank angle $\mu$ & $[0, 0.02]$ rad/s & 0.05 rad/s & $[-0.8, 0.8]$ rad \\
\midrule
gains & $k_v \in [0.2, 1]$ N/m & $k_\gamma \in [0.01, 0.04]$ & $k_\alpha \in [0.5, 2]$ 1/s, $k_\mu \in [0.05, 0.5]$ 1/s \\
length scales $\ell$ & references $[5, 30]$ s & noise $[2, 15]$ s & \\
scenario box & $z_0 \in [-6, 6]$ m & $v_0 \in [25, 30]$ m/s & $\gamma_0 \in [-0.06, 0.06]$ rad, coefficients $1.3 \times$ benchmark \\
flights & 60{,}000 & $t_f = 50$ s & $\Delta t = 0.1$ s \\
\bottomrule
\end{tabular}
\end{table}

The flight starts at the start state of the scenario with the initial inputs $(T_0, \alpha_0, \mu_0) = u_0$ of Table~\ref{tab:setup}. With $t_k = k\,\Delta t$, the state $(z_k, v_k, \gamma_k)$ at step $k$, and gains $k_v$, $k_\gamma$, $k_\alpha$ and $k_\mu$, the autopilot sets
\begin{equation}\label{eq:autopilot}
\begin{aligned}
  \alpha_{c,k} &= \alpha_{\mathrm{trim}}(v_k, z_k, \xi) + k_\gamma\,\big(\gamma_r(t_k) - \gamma_k\big), \\
  T_{k+1} &= T_k + \Delta t\, \big[\, k_v\,\big(v_r(t_k) - v_k\big) + n_T(t_k) \big]_{r_T}, \\
  \alpha_{k+1} &= \alpha_k + \Delta t\, \big[\, k_\alpha\,\big(\alpha_{c,k} - \alpha_k\big) + n_\alpha(t_k) \big]_{r_\alpha}, \\
  \mu_{k+1} &= \mu_k + \Delta t\, \big[\, k_\mu\,\big(\mu_r(t_k) - \mu_k\big) + n_\mu(t_k) \big]_{r_\mu},
\end{aligned}
\end{equation}
where $[\,\cdot\,]_r$ clips a rate to $[-r, r]$. The new inputs are clipped to their boxes, and the state advances by one RK4 step of \eqref{eq:uav} with the inputs linear within the step. The commanded angle of attack starts from the level-flight trim angle for the current speed, altitude and coefficients,
\begin{equation*}
  \alpha_{\mathrm{trim}} = -\frac{C_L^{\mathrm{req}} + C_{z0}}{C_{z\alpha}}, \qquad
  C_L^{\mathrm{req}} = \frac{2\, m_a g_0}{\varrho(z)\, v^2 S},
\end{equation*}
which follows from the small-angle lift $C_L \approx -(C_{z0} + C_{z\alpha}\,\alpha)$ in \eqref{eq:aero}. The speed and the flight-path angle of the generating flight are clipped to $[5, 60]$ m/s and $[-1.3, 1.3]$ rad, so that an unstable flight cannot overflow. The flight itself is discarded. Only the inputs at the 501 grid times and the scenario are stored.

The scenario of each flight is uniform in the training box of Table~\ref{tab:controls}, which is wider than the box of the benchmark by a factor of 1.3 about the centre for the aerodynamic coefficients. In training, each profile is paired with the scenario of its own flight or with a fresh draw from the training box, with probability one half each.

\subsection*{Training parameters}
The training loss is the residual \eqref{eq:residual} at 64 times per sample with causal weights of strength $\varepsilon = 0.1$ \citep{wang2024respecting}. We use the Adam optimizer for $6 \cdot 10^5$ steps with a batch size of 256. The learning rate rises linearly from 0 to $5 \cdot 10^{-4}$ over the first 2000 steps (warm-up) and then follows a cosine decay to $5 \cdot 10^{-6}$ at the last step. The training time was 2.4 h on one RTX 5070. The training parameters are summarized in Table~\ref{tab:training}.

\begin{table}[htbp]\centering\footnotesize
\caption{Training of the operator.}
\label{tab:training}
\begin{tabular}{cccccccc}
\toprule
training profiles & times per sample & $\varepsilon$ & batch & steps & step size & selected step & time \\
\midrule
60{,}000 generic & 64 & 0.1 & 256 & $6 \cdot 10^5$ & $5 \cdot 10^{-4}$ & $5.5 \cdot 10^5$ & 2.4 h \\
\bottomrule
\end{tabular}
\end{table}

The model selection is data-free. Among the snapshots of the training run, we keep the one with the smallest error estimate $\hat e$ of Section~\ref{sec:correction} on held-out generic controls, which is the snapshot after $5.5 \cdot 10^5$ steps. In the solver, the trunk values and their time derivatives on the correction grid are stored once before the optimization, and the correction of all $M$ scenarios runs by the parallel prefix scan of Section~\ref{sec:correction} in every iteration.

\subsection*{Initial guess and hardware}

All solvers start from the same initial guess and share the seeds 0, 2 and 3 of the sampled scenarios. The initial guess is the least-squares spline fit of the inputs that the autopilot \eqref{eq:autopilot} produces on the mean scenario with constant references of 27.5 m/s, zero flight-path angle and zero bank angle. This is a trim flight with nearly constant inputs. It holds no information about the target and has an RMS miss of about 2 km. The Adam solves run on one RTX 5070 and Ipopt on the 10 cores of the Intel Core Ultra 5 225 of workstations of the same type.

\subsection*{Reference Adam solver}
The reference Adam solver is the same single-shooting loop with the batched RK4 integration of \eqref{eq:uav} at step 0.5 s, 100 sequential steps per iteration for all $M$ scenarios, in place of the operator and the correction. The integration runs on the GPU.

\subsection*{Ipopt with exact Hessian}
Ipopt runs through CasADi with exact Hessians, the MUMPS linear solver, a tolerance of $10^{-5}$, and a filter line search. As an interior-point method, it handles the hard constraints with a logarithmic barrier whose weight it lowers monotonically toward zero, so the bounds hold at convergence. It also eliminates the states by RK4 integration at the same step of 0.5 s. It needs 55 to 65 iterations on the benchmark and 51 to 96 on the further tasks, and the build of the nonlinear program, 10 to 102 s at $M = 128$ to 4800, is not part of the solve time.

\subsection*{Ipopt with limited-memory Hessian}
Ipopt with a limited-memory Hessian in place of the exact one stops at the cap of 3000 iterations and is 11 to 12 times slower at $M = 128$ to 1600 with a slightly larger miss, so the exact Hessian is the fair and the fastest Ipopt configuration for these 75 variables.

\clearpage
\section{Reference solver and correction grid at other steps}\label{app:rk4steps}

Table~\ref{tab:rk4steps} repeats the benchmark solve at $M = 400$ with the reference Adam solver at RK4 steps of 0.1 to 2.5\,s and with the operator on correction grids of 0.5 and 1\,s, all with the margin of 0.035\,rad, as in \eqref{eq:margin}. The steps of 0.1 and 0.5\,s keep the mean flight-path angle below $\gamma_{\max}$ and give the same RMS miss as the operator. With steps of 1.0 and 2.5\,s, the solutions are infeasible, since the mean flight-path angle exceeds $\gamma_{\max}$ under RK4 at step 0.1\,s. The model gap of 2.8 to 195\,m shows that the solver underestimates its RMS miss on its own scenarios. The coarsest feasible step, 0.5\,s, takes 7.6 times the solve time of the operator.

The last row applies the correction of Section~\ref{sec:correction} to the reference Adam solver at step 1.0\,s, with the cubic Hermite interpolant of the RK4 solution as the prediction. The interpolant has zero residual at the RK4 nodes, so the correction needs a grid of 0.25\,s with four points per step. The corrected solver is then feasible and reaches the RMS miss of the reference solver at step 0.5\,s within 0.3\,m, at 4.2 times the solve time of the operator at $M = 400$, 3.2 times at $M = 1600$ and 2.6 times at $M = 4800$, where it is also 1.7 times slower than the reference solver at step 0.5\,s. At step 2.5\,s the prediction is outside the region in which the Newton step converges, and the corrected solver is exploited like the uncorrected one.

\begin{table}[htbp]\centering\footnotesize
\caption{Benchmark at $M = 400$, columns as in Table~\ref{tab:results}. A dagger marks a mean flight-path angle above $\gamma_{\max}$. The factor is relative to the operator on the 1\,s grid.}\label{tab:rk4steps}
\setlength{\tabcolsep}{4pt}
\begin{tabular}{@{}lcccccc@{}}
\toprule
solver & step or grid [s] & RMS miss [m] & model gap [m] & $\max_t \mathbb{E}[\gamma]$ & solve [s] & factor \\
\midrule
reference Adam & 0.1 & 52.6/52.7/52.5 & $<0.1$ & 0.492 & 817 & 38.6$\times$ \\
reference Adam & 0.5 & 52.6/52.6/52.4 & 0.1 & 0.512 & 161 & 7.6$\times$ \\
reference Adam & 1.0 & 55.0/55.1/55.8 & 2.8 & 0.683$^\dagger$ & 65 & 3.0$\times$ \\
reference Adam & 2.5 & 189.6/241.7/185.9 & 195.2 & 1.118$^\dagger$ & 38 & 1.8$\times$ \\
\midrule
operator & 0.5 & 52.7/52.7/52.3 & 0.1 & 0.492 & 32 & 1.5$\times$ \\
operator & 1.0 & 52.5/52.4/52.3 & $<0.1$ & 0.509 & 21 & 1.0$\times$ \\
\midrule
reference Adam, corrected & 1.0, grid 0.25 & 52.7/52.6/52.7 & 0.4 & 0.499 & 90 & 4.2$\times$ \\
\bottomrule
\end{tabular}

\end{table}

\clearpage
\section{Further experiments}\label{app:further}

\subsection*{Cost per iteration}
Table~\ref{tab:cost_iter} gives the cost of one iteration of the operator and of the reference Adam solver on the RTX 5070. The branch is one network evaluation at every $M$. The tensor contraction of the stored trunk values with the branch output is the linear term, bound by memory bandwidth, and the residual, the Jacobian and the scan are the correction. Above $M = 4800$, both costs grow at about 1.8 $\mu$s per scenario, so the factor between them falls slowly. The stored features take 313 KB per scenario, and at $M = 12{,}800$ the operator no longer fits the 12 GB of the GPU, while the reference solver runs.

\begin{table}[htbp]\centering\footnotesize
\caption{Cost of one iteration on the RTX 5070, profiled over 100 calls, and the stored trunk features of the operator.}\label{tab:cost_iter}
\setlength{\tabcolsep}{4pt}
\begin{tabular}{r*{7}{c}}
\toprule
 & \multicolumn{6}{c}{operator with the correction} & reference RK4, Adam \\
\cmidrule(lr){2-7}\cmidrule(lr){8-8}
$M$ & \shortstack{branch\\{}[ms]} & \shortstack{contraction\\{}[ms]} & \shortstack{residual and\\Jacobian [ms]} & \shortstack{scan\\{}[ms]} & \shortstack{value and\\gradient [ms]} & \shortstack{features\\{}[MB]} & \shortstack{value and\\gradient [ms]} \\
\midrule
400  & 0.04 & 0.21 & 0.03 & 0.06 & 1.0  & 125  & 8.3  \\
1600 & 0.04 & 0.85 & 0.08 & 0.18 & 3.0  & 501  & 10.5 \\
4800 & 0.03 & 2.58 & 0.27 & 0.83 & 8.0  & 1504 & 12.7 \\
6400 & 0.03 & 3.43 & 0.29 & 1.26 & 10.9 & 2005 & 15.3 \\
9600 & 0.03 & 5.19 & 0.50 & 1.94 & 16.6 & 3008 & 21.0 \\
\bottomrule
\end{tabular}
\end{table}

\subsection*{Time to accuracy}
Table~\ref{tab:time_to_accuracy} gives the held-out miss after a shorter budget of Adam steps. Both Adam solvers are then faster but worse and converge alike in steps, while the operator is 8 times faster at every level. The 2.4 h of training equals about 60 solves at $M = 400$ against the reference Adam solver and 20 against Ipopt, so the training pays off against Ipopt within Section~\ref{sec:exp} and against the Adam solver when the operator is reused.

\begin{table}[htbp]\centering\footnotesize
\caption{Time to accuracy at $M = 400$, means over three seeds. RMS miss on 200{,}000 held-out scenarios after a fixed number of Adam steps, and Ipopt at convergence.}\label{tab:time_to_accuracy}
\begin{tabular}{rcccc}
\toprule
 & \multicolumn{2}{c}{operator with the correction} & \multicolumn{2}{c}{reference RK4, Adam} \\
\cmidrule(lr){2-3}\cmidrule(lr){4-5}
Adam steps & RMS miss [m] & solve [s] & RMS miss [m] & solve [s] \\
\midrule
2000   & 68.5 & 2.1  & 66.7 & 15.9  \\
4000   & 54.9 & 4.2  & 55.0 & 31.8  \\
8000   & 53.4 & 8.5  & 53.6 & 63.4  \\
20{,}000 & 52.4 & 21.2 & 52.5 & 161.4 \\
\midrule
Ipopt  & \multicolumn{4}{c}{51.7 m in 413 s} \\
\bottomrule
\end{tabular}
\end{table}

\subsection*{Two correction steps for the operator}
A second correction step at the returned controls changes the RMS miss by less than 0.8 m and costs about 8 s more per solve. It lowers the RMS error of the corrected end state on the second target from between 4.1 and 5.9 m to 2.5 m and the largest error on the 40 s horizon from 26 to 4 m, while on the benchmark the grid and the stored features of the solver set a floor of 2 m that a second step does not lower.

\subsection*{Scenarios outside the training box}
On scenario boxes scaled by 1.3 (the training box) and 2 (outside it) about their centre, the operator matches the reference Adam solver within 0.3 and 1 m (the latter on the two seeds where that solver converged), the model gap grows from 0.1 to 1 m, and the largest value of the residual over the scenarios rises 3 to 4 times, so the residual flags scenarios outside the training range.

\end{document}